\documentclass[12pt]{amsart}
\usepackage[latin1]{inputenc}
\usepackage[english]{babel}
\usepackage{amsmath,amssymb,amsthm,amsfonts, color}
\usepackage{latexsym,dsfont}
\usepackage{enumerate}
\usepackage{graphicx}
\usepackage{latexsym}

\makeatletter
\@namedef{subjclassname@2010}{ \textup{2010} Mathematics Subject Classification}
\makeatother

\newtheorem{theorem}{Theorem}%[section]

\newtheorem{lemma}{Lemma}%[section]
\newtheorem{corollary}{Corollary}

\theoremstyle{definition}

\newtheorem{example}{Example}

\newtheorem{remark}{Remark}

\numberwithin{equation}{section}

\numberwithin{equation}{section}

\DeclareMathOperator*{\esssup}{ess\,sup}

\begin{document}
\title{Interpolation of abstract Ces\`aro, Copson \newline and Tandori spaces{\rm *}}
\thanks{{\rm *}This publication has been produced during scholarship period of the first author 
at the Lule{\aa} University of Technology, thanks to a Swedish Institute scholarschip (number 0095/2013).}

\author[Le\'snik]{Karol Le\'snik}
\address[Karol Le{\'s}nik]{Institute of Mathematics\\
of Electric Faculty, Pozna\'n University of Technology, ul. Piotrowo 3a, 60-965 Pozna{\'n}, Poland}
\email{\texttt{klesnik@vp.pl}}
\author[Maligranda]{Lech Maligranda}
\address[Lech Maligranda]{Department of Engineering Sciences and Mathematics\\
Lule{\aa} University of Technology, SE-971 87 Lule{\aa}, Sweden}
\email{\texttt{lech.maligranda@ltu.se}}

\begin{abstract}
We study real and complex interpolation of abstract Ces\`aro, Copson and Tandori spaces, including the description of 
Calder\'on-Lozanovski{\v \i} construction for those spaces. The results may be regarded as generalizations of interpolation 
for Ces\`aro spaces $Ces_p(I)$ in the case of real method, but they are new even for $Ces_p(I)$ in the case of complex 
method. Some results for more general interpolation functors are also presented. The investigations show an interesting 
phenomenon that there is a big difference between interpolation of Ces\`aro function spaces in the cases of finite and 
infinite interval.
\end{abstract}


\maketitle

%%%%%%%%%%%%%%%%%%%%%%%%%%%%%
\section{\protect \medskip Definitions and basic facts}
We recall some notations and definitions which will be needed. Let $(I,\Sigma,m)$ be a $\sigma$-complete measure space. By $L^0 = L^0(I)$ we denote the set of all equivalence 
classes of real-valued $m$ - measurable functions defined $I$. A {\it Banach ideal space} 
$X = (X, \|\cdot\|)$ is understood to be a Banach space contained in $L^0$, which satisfies the so-called ideal 
property: if $f, g \in L^0, |f| \leq |g|$ $m$-a.e. on $I$ and $g \in X$, then $ f\in X$ and $\|f\| \leq \|g\|$. Sometimes we write 
$\|\cdot\|_{X}$ to be sure which norm is taken in the space. If it is not stated otherwise we understand that in a Banach ideal 
space there is $f\in X$ with $f(x) > 0$ for each $x \in I$ (such an element is called the {\it weak unit} in $X$), which means that 
${\rm supp}X =I$. In the paper we concentrate on three underlying measure spaces $(I,\Sigma,m)$. If we say that $X$ is a 
{\it Banach function space} it means that it is a Banach ideal space where $I=[0,1]$ or $I=[0,\infty)$ and $m$ is just the Lebesgue  
measure, and $X$ is a {\it Banach sequence space} when  $I = \mathbb{N}$ with counting measure.
Later on by saying Banach ideal space we mean only one of those three cases.

For two Banach spaces $X$ and $Y$ the symbol $X\overset{A}{\hookrightarrow }Y$ means that the embedding $X \subset Y$ 
is continuous with the norm at most $A$, i.e., $\| f\|_{Y} \leq A \|f\|_{X}$ for all $f\in X$. When $X\overset{A}{\hookrightarrow }Y$ 
holds with some constant $A > 0$ we simply write $X\hookrightarrow Y$. Furthermore, $X = Y$ (or $X \equiv Y$) means that 
the spaces are the same and the norms are equivalent (or equal).

For a Banach ideal space $X = (X, \|\cdot\|)$ the {\it K{\"o}the dual space} (or {\it associated space}) $X^{\prime}$ is the 
space of all $f \in L^0$ such that the {\it associated norm}
\begin{equation} \label{dual}
\|f\|^{\prime} = \sup_{g \in X, \, \|g\|_{X} \leq 1} \int_{I} |f g | \, dm
\end{equation}
\vspace{-2mm}
is finite. The K{\"o}the dual $X^{\prime} = (X^{\prime}, \|\cdot \|^{\prime})$ is then a Banach ideal space.
Moreover, $X \overset{1}{\hookrightarrow }X^{\prime \prime}$ and we have equality $X = X^{\prime \prime}$ with 
$\|f\| = \|f\|^{\prime \prime}$ if and only if the norm in $X$ has the {\it Fatou property}, that is, if the conditions 
$0 \leq f_{n} \nearrow f$ a.e. on $I$ and $\sup_{n \in {\bf N}} \|f_{n}\| < \infty$ imply that $f \in X$ and $\|f_{n}\| \nearrow \|f\|$.

For a Banach ideal space $X = (X, \| \cdot\|)$ on $I$ with the K\"othe dual $X^{\prime}$ the following {\it generalized H\"older-Rogers
inequality} holds: if $f \in X$ and $g \in X^{\prime}$, then $fg$ is integrable and
\begin{equation} \label{1.2}
\int_I |f(x) g(x)| \, dx \leq \| f\|_X \| g \|_{X^{\prime}}.
\end{equation}

A function $f$ in a Banach ideal space $X$ on $I$ is said to have an {\it order continuous norm} in $X$ if, for any
decreasing sequence of $m$-measurable sets $A_{n} \subset I $ with $m(\bigcap A_n) = 0$, we have that  $\|f \chi_{A_{n}} \| \rightarrow 0$ 
as $n \rightarrow \infty$. The set of all functions in  $X$ with an order continuous norm is denoted by $X_{a}$. If $X_{a} = X$, then 
the space $X$ is said to be {\it order continuous} (we write shortly $X\in (OC)$). For an order continuous Banach ideal space $X$ 
the K{\"o}the dual $X^{\prime}$ and the dual space $X^{*}$ coincide. Moreover, a Banach ideal space $X$ with the Fatou property 
is reflexive if and only if both $X$ and its associate space $X^{\prime}$ are order continuous.

For a given {\it weight} $w$, i.e. a measurable function on $I$ with $0 < w(x) < \infty$ a.e. and for a Banach ideal space $X$ 
on $I$, the {\it weighted Banach ideal space} $X(w)$ is defined as $X(w)=\{f\in L^0: fw \in X\}$ with the norm 
$\|f\|_{X(w)}=\| f w \|_{X}$. Of course, $X(w)$ is also a Banach ideal space and 
\begin{equation} \label{dualweight}
[X(w)]^{\prime} \equiv X^{\prime}\Big(\frac{1}{w}\Big).
\end{equation}

By a {\it rearrangement invariant} or {\it symmetric space} on $I$ with the Lebesgue measure $m$, we mean a Banach 
function space $X=(X,\| \cdot \|_{X})$ with additional property that for any two equimeasurable functions 
$f \sim g, f, g \in L^{0}(I)$ (that is, they have the same distribution functions $d_{f}\equiv d_{g}$, where 
$d_{f}(\lambda) = m(\{x \in I: |f(x)|>\lambda \}),\lambda \geq 0)$, and $f\in E$ we have that $g\in E$ and $\| f\|_{E} = \| g\|_{E}$. 
In particular, $\| f\|_{X}=\| f^{\ast }\|_{X}$, where $f^{\ast }(t)=\mathrm{\inf } \{\lambda >0\colon \ d_{f}(\lambda ) < t\},\ t\geq 0$. 
Similarly one can define a {\it symmetric sequence space}. 
For general properties of Banach ideal spaces and symmetric spaces we refer to the books \cite{BS88}, \cite{KA77}, \cite{KPS82}, 
\cite{LT79} and \cite{Ma89}.

In order to define and formulate results we need the (continuous) Ces\`aro and Copson operators $C, C^*$ defined, respectively, 
as
$$
Cf(x) = \frac{1}{x} \int_0^x f(t) \,dt, 0 < x \in I, ~~C^*f(x) = \int_{I \cap [x, \infty)} \frac{f(t)}{t} \,dt,  x \in I,
$$
where $I=[0,1]$ or $I=[0,\infty)$.
The nonincreasing majorant $\widetilde{f}$ of a given function $f$,  is defined for $x \in I$ as
$$
\widetilde{f}(x) = \esssup_{t \in I, \, t \geq x} |f(t)|.
$$

For a Banach function space $X$ on $I$ we define the {\it abstract Ces\`aro (function) space} $CX = CX(I)$ as 
\begin{equation} \label{Cesaro}
CX=\{f\in L^0(I): C|f| \in X\} ~~ {\rm with ~the ~norm } ~~ \|f\|_{CX} = \| C|f| \|_{X},
\end{equation}
the {\it abstract Copson space} $C^*X = C^*X(I)$ as 
\begin{equation} \label{Copson}
C^*X=\{f\in L^0(I): C^*|f| \in X\} ~~ {\rm with ~the ~norm } ~~ \|f\|_{C^*X} = \| C^*|f| \|_{X},
\end{equation} 
and the  {\it abstract Tandori space} $\widetilde{X} = \widetilde{X} (I)$ as 
\begin{equation} \label{falka}
\widetilde{X}=\{f\in L^0(I): \widetilde{f}\in X\} ~~ {\rm with ~the ~norm } ~~ \|f\|_{\widetilde{X}}=\|\widetilde{f}\|_{X}.
\end{equation}

The {\it dilation operators} $\sigma_\tau$ ($\tau > 0$) defined on $L^0(I)$ by 
$$
\sigma_\tau f(x) = f(x/\tau) \chi_{I}(x/\tau) = f(x/\tau) \chi_{[0, \, \min(1, \, \tau)]}(x), ~~ x \in I,
$$
are bounded in any symmetric function space $X$ on $I$ and $\| \sigma_\tau \|_{X \rightarrow X} \leq \max \,(1, \tau)$ (see \cite[p. 148]{BS88} 
and \cite[pp. 96-98]{KPS82}). These operators are also bounded in some Banach function spaces which are not necessary symmetric. 
For example, if either $X = L^p(x^{\alpha})$ or $X = C(L^p(x^{\alpha}))$, then $\| \sigma_\tau\|_{X \rightarrow X} = \tau^{1/p + \alpha}$ 
(see \cite{Ru80} for more examples). 

In the sequence case the discrete Ces\`aro and Copson operators $C_d, C_d^*$ are defined for $n \in {\mathbb N}$ by 
$$
(C_d a)_n = \frac{1}{n} \sum_{k = 1}^n a_k, ~ (C^*_d a)_n = \sum_{k = n}^{\infty} \frac{a_k}{k}.
$$ 
The nonincreasing majorant $\widetilde{a}$ of a given sequence $a = (a_n)$ is defined for $n \in \mathbb N$ as
$$
\widetilde{a_n} = \sup_{k \in {\mathbb N}, \, k \geq n} |a_k|.
$$
Then the corresponding {\it Ces\`aro sequence space} $C_dX$, {\it Copson sequence space} $C_d^*X$ 
and {\it Tandori sequence space} $\widetilde{X}_d$ are defined analogously as in (\ref{Cesaro}), (\ref{Copson}) 
and (\ref{falka}). Moreover, for every $m \in \mathbb N$ let $\sigma_m$ and $\sigma_{1/m}$ be the {\it dilation operators} 
defined in spaces of sequences $a = (a_n)$ by (cf. \cite[p. 131]{LT79} and \cite[p. 165]{KPS82}):
$$
\sigma_m a = \left( ( \sigma_m a)_n \right)_{n=1}^{\infty} = \big (a_{[\frac{m-1+n}{m}]} \big)_{n=1}^{\infty} 
= \big ( \overbrace {a_1, a_1, \ldots, a_1}^{m}, \overbrace {a_2, a_2, \ldots, a_2}^{m}, \ldots \big)
$$
\begin{eqnarray*}
\sigma_{1/m} a 
&=& 
\left( ( \sigma_{1/m} a)_n \right)_{n=1}^{\infty} = \Big (\frac{1}{m} \sum_{k=(n-1)m + 1}^{nm} a_k \Big)_{n=1}^{\infty} \\
&=& 
\big ( \frac{1}{m} \sum_{k=1}^m a_k, \frac{1}{m} \sum_{k=m+1}^{2m} a_k, \ldots, \frac{1}{m} \sum_{k=(n-1)m + 1}^{nm} a_k, \ldots \big).
\end{eqnarray*}
These operators are discrete analogs of the dilation operators $\sigma_{\tau}$ defined in function spaces.
They are bounded in any symmetric sequence space but also in some Banach sequence spaces, for example, 
$\| \sigma_{m} \|_{l^p(n^{\alpha}) \rightarrow l^p(n^{\alpha})} \leq m^{1/p} \max(1, m^{\alpha})$ and 
$\| \sigma_{1/m} \|_{l^p(n^{\alpha}) \rightarrow l^p(n^{\alpha})} \leq m^{-1/p} \max \,(1, m^{-\alpha})$.

Properties of Ces\`aro sequence spaces $ces_p = Cl^p$ were investigated in many papers (see \cite{MPS07} and references 
given there), while properties of Ces\`aro function spaces $Ces_p(I) = CL^p(I)$ we can find in \cite{AM09} and \cite{AM14b}. 
Abstract Ces\`aro spaces $CX$ for Banach ideal spaces $X$ on $[0, \infty)$ were defined already in \cite{Ru80} and spaces 
$CX, \,\widetilde{X}$ for $X$ being a symmetric space on $[0,\infty)$ have appeared, for example, in \cite{KMS07}, \cite{DS07} 
and \cite{AM13}. General considerations of abstract Ces\`aro spaces began to be studied in papers \cite{LM15a, LM15b}.

Copson sequence spaces $cop_p = C^*l^p$ and Copson function spaces $Cop_p = C^*L^p$ on $[0, \infty)$ we can find in 
G. Bennett's memoir \cite[pp. 25-28 and 123]{Be96}. Moreover, Copson function spaces $Cop_p = C^*L^p$ on $[0, 1]$ were 
used in \cite{AM13} (see also \cite{AM14a}, \cite{AM14b}) to understand Ces\`aro spaces and their interpolation.

We put the name generalized Tandori function spaces on $\widetilde{X} = \widetilde{X}(I)$ in honour of Tandori who proved 
in 1954 that the dual space to $(Ces_{\infty}[0, 1])_a$ is $\widetilde{L^1[0, 1]}$. In 1966 Luxemburg-Zaanen \cite[Theorem 4.4]{LZ66} have
found the K\"othe dual of  $Ces_{\infty}[0, 1]$ as $(Ces_{\infty}[0, 1])^{\prime} \equiv \widetilde{L^1[0, 1]}$. Already in 1957,  
Alexiewicz \cite[Theorem 1]{Al57} showed (even for weighted case) that $( \widetilde{l^1} )^{\prime} \equiv ces_{\infty}$. 
In \cite[Theorem 7]{LM15a} we were able to give a simple proof of a generalization of the Luxemburg-Zaanen duality theorem: 
$\big [ C(L^{\infty}(v))  \big]^{\prime} \equiv \widetilde{L^1(w)}$, where $v(x) = x/\int_0^x w(t)\, dt, x \in I$. 
Bennett \cite[Corollary 12.17]{Be96} proved that $(ces_p)^{\prime} = \widetilde{l^{p^{\prime}}}$ for $1 < p < \infty$. Surprisingly, the dual 
of  Ces\`aro function space is essentaially different for $I = [0, \infty)$ and for $I = [0,1]$, as it was proved by Astashkin-Maligranda \cite{AM09} (see also \cite{LM15a} for a simpler proof). 
Namely, for $1 < p < \infty$ we have $(Ces_p[0, \infty))^{\prime} = \widetilde{L^{p^{\prime}}[0, \infty)}$ (cf. \cite[Theorem 2]{AM09}) and  
$(Ces_p[0, 1])^{\prime} = \widetilde{L^{p^{\prime}}(\frac{1}{1-x})[0, 1]}$ (cf. \cite[Theorem 3]{AM09}). Generalized Tandori 
spaces $\widetilde{X}$ (without this name on the spaces) and their properties appeared in papers  \cite[p. 935]{LM15a} 
and \cite[p. 228]{LM15b}.
\vspace{3mm}

We will need the following result on duality of abstract Ces\`aro spaces proved in \cite{LM15a}.
\vspace{1mm}

%%%%%%%%%%%%%%%%%%%%% Theorem A
{\bf Theorem A.} {\it 
Let $X$ be a Banach ideal space with the Fatou property such that the Ces\`aro operator $C$ is bounded on $X$.
\begin{itemize}
\item[(i)] If $I = [0,\infty)$ and the dilation operator $\sigma_\tau$ is bounded on $X$ for some $0 < \tau < 1$, then 
\begin{equation} \label{Thm2i}
(CX)^{\prime} = \widetilde{X^{\prime}}.
\end{equation}
\item[(ii)] If $X$ is a symmetric space on $[0, 1]$ such that $C, C^*: X \rightarrow X$ are bounded, then
\begin{equation} \label{Thm2ii}
(CX)^{\prime} = \widetilde{ X^{\prime}(1/v)}, ~ {\rm where} ~~ v(x) = 1-x, ~x \in [0, 1).
\end{equation}
\item[(iii)] If $X$ is a sequence space and the dilation operator $\sigma_{3}$ is bounded on $X^{\prime}$, then 
\begin{equation} \label{Thm2iii}
(CX)^{\prime} = \widetilde{X^{\prime}}.
\end{equation}
\end{itemize}
}
%%%%%%%%%%%%%%%%%%%%%%%%%
\vspace{3mm}

The paper is organized as follows. In Section 2 we present comparisons of Ces{\`a}ro, Copson and Tandori spaces 
as well as the ``iterated'' spaces $CCX$ and $C^*C^*X$. 

Section 3 contains results on commutativity of the Calder\'on-Lozanovski{\v \i} construction with abstract Ces{\`a}ro spaces 
and with generalized Tandori spaces. There are some differences in assumptions on Ces{\`a}ro function spaces in the 
cases on $[0, \infty)$ and on $[0, 1]$, and the sequence case, as we can see in Theorem \ref{thm:cesCL}.
Important in our investigations were results proved in \cite{LM15a} on  K{\"o}the duality of $CX$ (cf. Theorem A). 
In the case of generalized Tandori spaces we were able to prove an analogous result in Theorem \ref{thm:tyldaspace} by using another method, that ommits the duality argument.

Results proved here are then used to described interpolation of abstract Ces{\`a}ro and Tandori spaces by the complex method. 
Identifications in Theorem 5 are new even for classical Ces{\`a}ro spaces $Ces_p(I)$.

In Section 4, the commutativity of the real method of interpolation with abstract Ces{\`a}ro spaces is investigated in 
Theorem \ref{thm:Cesaro=real}. We also collected here our knowledge about earlier results on interpolation of Ces{\`a}ro 
spaces $Ces_p$ and their weighted versions. 

Finally, in Section 5, we collected information on Calder\'on couples of Ces{\`a}ro spaces and some related spaces. Several remarks 
and open problems are also formulated. From all the above discussions we can see a big difference between interpolation of abstract 
Ces{\`a}ro spaces on intervals $[0, \infty)$ and $[0, 1]$.

%%%%%%%%%%%%%%%%%%%%%%%%% Section 2
\section{Comparison of Ces{\`a}ro, Copson and Tandori spaces}

First of all notice that  $X\overset{A}{\hookrightarrow }Y$ implies  $CX\overset{A}{\hookrightarrow }CY$, 
$C^*X\overset{A}{\hookrightarrow }C^*Y$ and $\widetilde{X} \overset{A}{\hookrightarrow } \widetilde{Y}$. 
Moreover, it can happend that spaces are different but corresponding Ces{\`a}ro, Copson and Tandori spaces 
are the same, that is, there are $X \neq Y$ such that $CX = CY$, $C^*X = C^*Y$ and $ \widetilde{X} =  \widetilde{Y}$.

%%%%%%%%%%%%% Example 1
\begin{example} \label{Ex1} 
If $X = L^2[0, \frac{1}{4}] \oplus L^{\infty}[\frac{1}{4}, \frac{1}{2}] \oplus L^2[\frac{1}{2}, 1]$, then $X\hookrightarrow L^2[0, 1]$,  
$CX = CL^2 = Ces_2[0,1], C^*X = C^*L^2 = Cop_2[0, 1]$ and $ \widetilde{X} =  \widetilde{L^2[0, 1]}$, because 
\begin{eqnarray*}
\sup_{x \in [ \frac{1}{4},  \frac{1}{2}]} \frac{1}{x} \int_0^x |f(t)| \, dt 
&\leq &
4 \int_0^{1/2} |f(t)| \, dt = 4 \int_0^{1/2} |f(t)| \, dt  \, (\int_{1/2}^1 x^{-2} dx)^{1/2} \\
&\leq&
4 \left[ \int_{1/2}^1 \big( \frac{1}{x} \int_0^x |f(t)| \, dt \big)^2dx \right]^{1/2},
\end{eqnarray*}
\begin{equation*}
\sup_{x \in [ \frac{1}{4},  \frac{1}{2}]} \int_x^1 \frac{|f(t)|}{t} \, dt =  \int_{1/4}^1 \frac{|f(t)|}{t} \, dt \leq 
2 \left[ \int_0^{1/4} \big( \int_x^1 \frac{|f(t)|}{t} \, dt \big)^2 dx\right]^{1/2},
\end{equation*}
and
\begin{equation*}
\sup_{x \in  [ \frac{1}{4},  \frac{1}{2}]}  \widetilde{f}(x) = \widetilde{f}(1/4) \leq 2 \, \big( \int_0^{1/4}  \widetilde{f}(x)^2 \, dx \big)^{1/2}.
\end{equation*}
\end{example}
%%%%%%%%%%%%%

Since $\widetilde{X} \overset{1}{\hookrightarrow } X$ it follows that $C \widetilde{X} \overset{1}{\hookrightarrow }CX$ and the reverse 
imbedding, under some assumptions on $X$, was proved in \cite[Theorem 1]{LM15b}. Namely, consider the maximal operator $M$ 
(defined for $x \in I$ by $Mf(x) = \sup_{a, b \in I, 0 \leq a \leq x \leq b} \frac{1}{b-a} \int_a^b |f(t)|\, dt$) and a Banach ideal space $X$ 
on $I$. In the case $I = [0, \infty)$, if $M$ is bounded on $X$, then
\begin{equation} \label{2.1}
CX \overset{B}{\hookrightarrow } C \widetilde{X} ~~ {\rm with} ~ B = 4 \, \| M \|_{X \rightarrow X} 
\end{equation}
(cf. \cite[Theorem 1(i)]{LM15b}), and in the case $I = [0, 1]$ if $M, \sigma_{1/2}$ are bounded on $X$ and $L^{\infty} \hookrightarrow X$, then
\begin{equation} \label{2.2}
CX \cap L^1 \overset{B_1}{\hookrightarrow } C \widetilde{X}  \overset{B_2}{\hookrightarrow } CX \cap L^1
\end{equation}
with $B_1 = 4 \, \| M \|_{X \rightarrow X}  \| \sigma_{1/2} \|_{X \rightarrow X}$ and $B_2 = \max \,(1, 1/\| \chi_{[0, 1]} \|_X)$ (cf. 
\cite[Theorem 1(iii)]{LM15b}). In particular, if $X$ is a symmetric space on $I$ and $C$ is bounded on $X$, then
\begin{equation} \label{2.3}
C \widetilde{X} = CX  ~ {\rm for} ~ I = [0, \infty) ~ {\rm and} ~~ C \widetilde{X} = CX \cap L^1 ~ {\rm for} ~ I = [0, 1].
\end{equation}
We should mention here that the boundedness of $C$ on a symmetric space $X$ implies boundedness of the maximal operator 
$M$ on $X$, which follows from the Riesz inequality $(Mf)^*(x) \leq c\, Cf^*(x)$ true for any $x \in I$ with a constant $c \geq 1$ independent 
of $f$ and $x$ (cf. \cite[p. 122]{BS88}).
\vspace{0.1mm}

Now, we collect inclusions and equalities between Ces{\`a}ro spaces $CX$, Copson spaces $C^*X$, their iterations $ CCX, C^*C^*X$ 
and Tandori spaces $ \widetilde{X}$. Some results were proved before for $X = l^p$ by Bennett in \cite{Be96} and for $X = L^p$ by 
Astashkin-Maligranda \cite{AM09}. Moreover,  Curbera and Ricker in  \cite{CR13} already proved point (viii) in the theorem below. 
Let us recall that the {\it unilateral shift} $S$ on a sequence space is defined by $S(x_1, x_2, x_3, \ldots) = (0, x_1, x_2, x_3, \ldots)$.

%%%%%%%%%%%%%%%%%Theorem 1
\begin{theorem}\label{thm:relations} \item[(a)] Let $X$ be a Banach function space on $I = [0, \infty)$.
\begin{itemize}
\item[(i)] If $C$ is bounded on $X$, then $X \hookrightarrow CX$ and if, additionally, the dilation operator 
$\sigma_{1/a}$ is bounded on $X$ for some $a > 1$, then $CX = CCX$.
\item[(ii)] If $C^*$ is bounded on $X$, then $X \hookrightarrow C^*X$ and if, additionally, the dilation operator 
$\sigma_{1/a}$ is bounded on $X$ for some $0 < a < 1$, then $C^*X = C^*C^*X$.
\item[(iii)] If both operators $C$ and $C^*$ are bounded on $X$, then $CX = C^*X$. 
\end{itemize}
\item[(b)]  Let $X$ be a Banach function space on $I = [0, 1]$.
\vspace{-2mm}
\begin{itemize}
\item[(iv)] If $C$ is bounded on $X$, then $X \overset{A}{\hookrightarrow } CX$ and $C^*X\overset{A}{\hookrightarrow } CX$ 
with $A = \| C \|_{X \rightarrow X}$. The last embedding is strict even if $X = L^p[0, 1]$ with $1 < p < \infty$. 
\item[(v)] If $C^*$ is bounded on $X$, then $X \hookrightarrow C^*X$. If, additionally, the dilation operator $\sigma_{1/a}$ is
bounded on $X$ for some $0 < a < 1$, then $C^*X = C^*C^*X$.
\item[(vi)] If $C^*$ is bounded on $X$ and $L^{\infty} \hookrightarrow X$, then $C^*X \hookrightarrow CX \cap L^1$. If, 
additionally, $C$ is bounded on $X$ and $X \hookrightarrow L^1$, then  $C^*X= CX \cap L^1$ and 
$(C^*X)^{\prime} = \widetilde{X^{\prime}}$.
\item[(vii)] If both operators $C$ and $C^*$ are bounded on $X$, and $X$ is a symmetric space, then 
$C^*X = CX \cap L^1 = C \widetilde{X}$.
\end{itemize}
\item[(c)] Let $X$ be a Banach sequence space.
\begin{itemize}
\item[(viii)] If $C_d$ is bounded on $X$, then $X \hookrightarrow C_dX$ and if, additionally, the dilation operator $\sigma_{1/2}$ 
is bounded on $X$, then $C_dX = C_dC_dX$.
\item[(ix)] If $C_d^*$ is bounded on $X$, then $X \hookrightarrow C_d^*X$ and if, additionally, the dilation operator 
$\sigma_2$ is bounded on $X$, then $C_d^*X = C_d^*C_d^*X$.
\item[(x)] If operators $C_d, C_d^*$ and unilateral shift $S$ with its dual $S^*$ are bounded on $X$, then $C_dX = C_d^*X$.
\end{itemize}
\end{theorem}
%%%%%%%%%%%%%%%%%%%%

%%%%  Proof
\proof
(i) The first inclusion is clear from which we obtain also $CX \hookrightarrow CCX$. Then the equality $CX = CCX$ 
follows from Lemma 6 in \cite{LM15a}, where it was proved that $CC|f|(x) \geq \frac{\ln a}{a} \, C|f|(x/a)$ for all $x > 0$ 
with arbitrary $a > 1$. Thus 
$$
\sigma_{1/a} CC|f|(x) = CC|f|(ax) \geq  \frac{\ln a}{a} \, C|f|(x),
$$
and so
$$
\| C|f|\|_X \leq  \frac{a}{\ln a} \| \sigma_{1/a} CC|f| \|_X \leq  \frac{a}{\ln a} \| \sigma_{1/a}\|_{X \rightarrow X} \| CC|f| \|_X,
$$
which gives the required inclusion $CCX \hookrightarrow CX$. 

(ii) Once again the first inclusion comes directly from the assumption and thus $C^*X \hookrightarrow C^*C^*X$. 
To get the reverse inclusion observe that for $f \geq 0, x>0$ and $0 < a < 1$ by monotonicity of $C^*f$ we have
\begin{eqnarray*}
\sigma_{1/a}C^*C^*f(x) 
&=& 
C^*C^*f(ax) = \int_{ax}^{\infty} \frac{C^*f(t)}{t} \, dt \\
&\geq& 
\int_{ax}^x\frac{C^*f(t)}{t} \, dt \geq C^*f(x) \ln \frac{1}{a}.
\end{eqnarray*}
Thus,
$$
\| C^* f\|_X \leq \frac{1}{\ln1/a} \| \sigma_{1/a}(C^*C^*)f \|_X \leq \frac{\| \sigma_{1/a} \|_{X \rightarrow X}}{\ln1/a} \, \| C^*C^*f \|_X,
$$
which gives the required inclusion $C^*C^*X \hookrightarrow C^*X$.

\vspace{2mm}

(iii) Since $C|f| + C^*|f| = C^*C|f|$ and $C|f| + C^*|f| = CC^*|f|$ it follows that
\begin{eqnarray*}
\| f \|_{CX} 
&=& 
\| C|f| \|_X \leq \| C|f| + C^*|f| \|_X = \| C^*C|f| \|_X \\
&\leq&
\| C^* \|_{X \rightarrow X} \| C|f| \|_X = \| C^* \|_{X \rightarrow X} \, \| f \|_{CX} 
\end{eqnarray*}
and
\begin{eqnarray*}
\| f \|_{C^*X} 
&=& 
\| C^*|f| \|_X \leq \| C|f| + C^*|f| \|_X = \| CC^*|f| \|_X \\
&\leq&
\| C \|_{X \rightarrow X} \| C^*|f| \|_X = \| C \|_{X \rightarrow X} \, \| f \|_{C^*X}. 
\end{eqnarray*}
Therefore,
\begin{eqnarray*}
\| f \|_{CX} 
&\leq& 
\| C|f| + C^*|f| \|_X \leq \| C \|_{X \rightarrow X} \, \| f \|_{C^*X} \\
&\leq&
 \| C \|_{X \rightarrow X} \, \| C|f| + C^*|f| \|_X \leq \| C \|_{X \rightarrow X} \, \| C^* \|_{X \rightarrow X} \, \| f \|_{CX}, 
\end{eqnarray*}
and so $C^*X\overset{A}{\hookrightarrow } CX \overset{ B}{\hookrightarrow } C^*X$ with $A = \| C \|_{X \rightarrow X}$ and 
$B = \| C^* \|_{X \rightarrow X}$. 

\vspace{2mm}

(iv) The first inclusion is clear. The second embedding follows from the equality $ C + C^* = C C^*$, which gives
\begin{equation} \label{imbeddingCC}
C^*X\overset{A}{\hookrightarrow } CX ~ {\rm with} ~ A = \| C \|_{X \rightarrow X}.
\end{equation}
Moreover, equality in (\ref{imbeddingCC}) does not hold in general as it was shown already in \cite[p. 48]{AM09} for $X = L^p[0, 1]$.
\vspace{2mm}

(v) The prove is the same as in $(ii)$.
\vspace{2mm}

(vi) Since for $f \geq 0$ and $x \in [0, 1]$ we have (see also \cite[p. 48]{AM13}) 
$$
C^*Cf(x) = Cf (x) + C^*f(x) - \int_0^1 f(t) \, dt
$$
it follows that
\begin{eqnarray*}
\| f \|_{C^*X} 
&=& 
\| C^*f \|_X \leq \| C^*f + Cf \|_X =  \| C^* Cf + \int_0^1 f(t) \, dt \, \chi_{[0, 1]} \|_X \\
& \leq& 
\| C^* \|_{X \rightarrow X} \| Cf\|_X + \| f \|_{L^1}  \, \| \chi_{[0, 1]} \|_X  \\
&\leq&
4 \, \max \{\| C^* \|_{X \rightarrow X},  \| \chi_{[0, 1]} \|_X \}  \, \max \{ \| Cf\|_X, \| f \|_{L^1}  \}\\
&=&
4 \, \max \{\| C^* \|_{X \rightarrow X},  \| \chi_{[0, 1]} \|_X \}  \, \| f\|_{CX \cap L^1},
\end{eqnarray*}
that is, $CX \cap L^1\overset{D}{\hookrightarrow } C^*X$ with $D = 4 \, \max \{\| C^* \|_{X \rightarrow X},  \| \chi_{[0, 1]} \|_X \}$.

On the other hand, from (\ref{imbeddingCC}) and since $X \hookrightarrow L^1$ it follows $C^*X \hookrightarrow C^*L^1 \equiv L^1$ 
and we obtain $C^*X \hookrightarrow CX \cap L^1$. Therefore, $C^*X = C^*L^1 \cap L^1$. 

The embedding $(C^*X)^{\prime} \hookrightarrow  \widetilde{X^{\prime}}$ can be proved in the following way. 
Using just mentioned identification (\ref{2.3}), equality of the K\"othe dual of the sum as the intersection of K\"othe duals 
(see, for example, \cite[Lemma 3.4, p. 342]{LZ66} or \cite[Lemma 15.5]{Ma89}) and Theorem A(ii) we obtain
$$
(C^*X)^{\prime} = (CX \cap L^1)^{\prime} = (CX)^{\prime} + L^{\infty} = \widetilde{X^{\prime}(1/v)} + L^{\infty}.
$$
Then, since $v(x) = 1 - x \leq 1$ it follows that $\widetilde{X^{\prime}(1/v)} \overset{1}{\hookrightarrow }  \widetilde{X'}$ and by the assumption 
$X \hookrightarrow L^1$ we obtain $L^{\infty} \hookrightarrow  X^{\prime}$ which gives 
$L^{\infty} = \widetilde{L^{\infty}}  \hookrightarrow  \widetilde{X^{\prime}}$. Thus, $\widetilde{X^{\prime}(1/v)} + L^{\infty} \hookrightarrow  \widetilde{X'}$.

To finish the proof we need to show the embedding $\widetilde{X'} \hookrightarrow \widetilde{X^{\prime}(1/v)} + L^{\infty}$. Let $0\leq f\in \widetilde{X'}$. 
Then
$$
\|(f - \tilde f(1/2))_+\|_{\widetilde{X^{\prime}(1/v)}}=\|\frac{(\tilde f - \tilde f(1/2))_+}{v}\|_{X^{\prime}} 
\leq \frac{1}{1-1/2}\|(\tilde f - \tilde f(1/2))_+\|_{X^{\prime}}\leq  2 \, \|\tilde f \|_{X^{\prime}}.
$$
Moreover, by the H\"older-Rogers inequality (\ref{1.2}), we obtain for $f \in X$ and any $0 < t < 1$,
$$
f^*(t) \leq \frac{1}{t} \int_0^t f^*(s)\, ds \leq \frac{1}{t} \| \chi_{[0, t]}\|_X \, \| f^*\|_{X^{\prime}}.
$$
Therefore, $\|\tilde f(1/2) \chi_{[0,1]}\|_{L^{\infty}} = \tilde f(1/2) \leq 2 \, \varphi_X(1/2) \|\tilde f\|_{X'} $,
and consequently 
$$
\| f \|_{\widetilde{X^{\prime}(1/v)} + L^{\infty}} \leq [2+ 2\, \varphi_{X}(1/2)] \, \| f \|_{\widetilde {X'}} \,.
$$

(vii) The first equality follows from (vi) and the second from (\ref{2.3}).

(viii) Of course,  $C_dX \hookrightarrow C_d C_dX$ and the reverse inclusion for $X = l^p$ was already given by 
Bennett \cite[20.31]{Be96}, but it was simplified and generalized by Curbera and Ricker \cite[Proposition 2]{CR13} 
who have shown that for $n \geq 2$ there holds
\begin{equation} \label{CR}
\dfrac{1}{[n/2]} \sum_{j=1}^{[n/2]} |a_j| \leq 6 \sum_{k=1}^n \dfrac{1}{k}  \sum_{j=1}^k |a_j|.
\end{equation}
Thus $(C_d a)_n \leq 6\, (C_d C_d a)_{2n} \leq 12 \, (\sigma_{1/2} C_d C_d a)_n$ and 
$$
\| C_d a \|_X \leq 12 \|\sigma_{1/2}\|_{X \rightarrow X} \, \| C_dC_d a\|_X,
$$
which gives the required inclusion.

(ix) Of course,  $C^*_dX \hookrightarrow C^*_d C^*_dX$ and we need only to prove the reverse inclusion. Since
$ \sum_{k=n}^{2n-1} \frac{1}{k} \geq \int_n^{2n} \frac{1}{t} dt = \ln 2 \geq \frac{1}{2}$ it follows that
\begin{eqnarray*}
(C^*_d C^*_d a)_n 
&=&
\sum_{k=n}^{\infty} \frac{(C^*_d a)_k}{k} \geq \sum_{k=n}^{2n-1} \frac{(C^*_d a)_k}{k} \geq
(C^*_d a)_{2n-1} \sum_{k=n}^{2n-1} \frac{1}{k} \\
&\geq&
\frac{1}{2} \, (C^*_d a)_{2n-1} \geq \frac{1}{2} \, (C^*_d a)_{2n},
\end{eqnarray*}
and
$$
(\sigma_2C^*_d C^*_d a)_n = (C^*_d C^*_d a)_{[\frac{n+1}{2}]} \geq \frac{1}{2} \, (C^*_d a)_n.  
$$
Thus,
$$
\| C^*_d a \|_X \leq 2 \, \| \sigma_2 \|_{X \rightarrow X} \, \| C^*_d C^*_d a\|_X,
$$
which gives the required inclusion.

(x) This result for $X = l^p$ was proved by Bennett \cite[p. 47]{Be96} who observed that
\begin{equation*}
C_d = (C_d - S^*) C^*_d ~~{\rm and} ~~ C^*_d = (C^*_d - I) S C_d,
\end{equation*}
where $S$ is the unilateral shift and $S^*$ its dual.
Of course, his proof is also working for more general Banach sequence spaces. Namely,
$$
\| C_d a \|_X \leq \left( \| C_d \|_{X \rightarrow X} + \| S^* \|_{X \rightarrow X} \right) \| C^*_d a \|_X = A \, \| C^*_d a \|_X,
$$
which gives $C_d^*X \overset{A}\hookrightarrow C_dX$. Also
$$
\| C_d^* a \|_X \leq \left( \| C_d S\|_{X \rightarrow X} + \| S \|_{X \rightarrow X} \right) \| C^*_d a \|_X = B \, \| C^*_d a \|_X
$$
and so  $C_dX \overset{B}\hookrightarrow C_d^*X$. Putting together both inclusions we obtain $C_dX = C_d^*X$.
\endproof
%%%%%%%%%%%%%%%%Remark 1
\begin{remark} \label{Re1}
On $[0,1]$ the space $CXX$ may be essentially bigger than $CX$. In fact, taking $f(x) = \frac{1}{(1 - x)^2}, x \in(0,1)$ we have 
$f\in CCL^p$ for any $1 \leq p < \infty$ and $f \not \in CL^p$ since $Cf(x) = \frac{1}{1-x} \not \in CL^1$. Moreover, if the  operator 
$C$ is not bounded in $X$ on $[0, 1]$, then embedding  relationships between $X, CX$ and $CCX$ may not hold. 
For $X = L^1[0, 1]$ we have
\begin{equation*}
\| f \|_{CL^1} = \int_0^1 |f(x)| \, \ln (1/x) \, dx ~ {\rm and} ~ \| f \|_{CCL^1} = \int_0^1 |f(x)| \, \ln^2 (1/x) \, dx.
\end{equation*}
Spaces $L^1, CL^1$ and $CCL^1$ are not comparable. In fact, for $0 < \alpha < 1$ let $f_{\alpha}(x) = \frac{1}{x} \chi_{[\alpha, 1]}$. 
Then $\| f_{\alpha} \|_{L^1} = \ln (1/\alpha), \| f_{\alpha} \|_{CL^1} = \frac{1}{2} \ln^2 (1/\alpha), 
\| f_{\alpha} \|_{CCL^1} = \frac{1}{3} \ln^3 (1/\alpha)$ and
$$
\dfrac{2 \, \| f_{\alpha} \|_{CL^1}}{\| f_{\alpha}  \|_{L^1}} = \dfrac{3 \, \| f_{\alpha} \|_{CCL^1}}{\| f_{\alpha} \|_{CL^1}} =  \ln \frac{1}{\alpha} \rightarrow \infty ~{\rm as} ~ \alpha \rightarrow 0^+ ~( {\rm and} ~ \rightarrow 0 ~ {\rm as} ~ \alpha \rightarrow 1^-).
$$
\end{remark}

%%%%%%%%%%%%%%%%%%%%%%%%% Section 2
\section{Calder\'on-Lozanovski{\v \i} construction}

Let us recall the Calder\'on-Lozanovski{\v \i} construction for Banach ideal spaces. The class ${\mathcal U}$ consists of all functions 
$\varphi: {\mathbb R_+} \times {\mathbb R_+} \rightarrow {\mathbb R_+}$ that are positively homogeneous (i.e., 
$\varphi(\lambda s, \lambda t) = \lambda \varphi (s, t)$ for every $s, t, \lambda \geq 0$) and concave, that is 
$\varphi (\alpha s_1 + \beta s_2, \alpha t_1 + \beta t_2) \geq \alpha \varphi (s_1, t_1) + \beta  \varphi(s_2, t_2)$ for all 
$\alpha, \beta \in [0, 1]$ with $\alpha + \beta = 1$, and all $s_i, t_i \geq 0, i = 0, 1$. Note that any function $\varphi \in {\mathcal U}$ is 
continuous on $(0, \infty) \times (0, \infty)$.

Given such $\varphi \in {\mathcal U}$ and a couple of Banach ideal spaces $(X_0, X_1)$ on the same measure space, the 
{\it Calder\'{o}n-Lozanovski{\u \i} space} $\varphi (X_0, X_1)$ is defined as the set of all $f \in L^0$ such that for some $f_0 \in X_0, f_1 \in X_1$ 
with $\| f_0 \|_{X_0} \leq 1, \| f_1\|_{X_1} \leq 1$ and for some $\lambda > 0$ we have
\begin{equation*}
| f(x)| \leq \lambda\, \varphi(|f_0(x)|, | f_1(x)|) ~~ {\rm a.e. ~on} ~ I. 
\end{equation*}
The norm $ \| f \|_{\varphi (X_0, X_1)}$ of an element $f \in \varphi(X_0, X_1)$ is defined as the 
infimum of those values of $\lambda$ for which the above inequality holds and the space $(\varphi (X_0, X_1), \| \cdot \|_{\varphi})$ 
is then a Banach ideal space. It can be shown that
\begin{equation*}
\varphi (X_0, X_1) = \left\{ f\in L^{0}: |f| \leq \varphi (f_0, f_1)  ~{\rm for ~some} ~ 0 \leq f_0 \in X_0, 
~ 0 \leq f_1 \in X_1 \right \}
\end{equation*}
with the norm 
\begin{equation*}
\| f \|_{\varphi \left( X_0, X_1\right) }=\inf \left\{ \max \left\{\| f_0 \|_{X_0}, \| f_1 \|_{X_1}\right\}: \text{\ }
| f| \leq \varphi (f_0, f_1)\text{,}\ 0 \leq f_0 \in X_0, 0 \leq f_1 \in X_1 \right\} \text{.}  
\end{equation*}
In the case of power functions $\varphi (s, t) = s^{\theta } t^{1-\theta} $ with $0 \leq \theta \leq 1$ spaces $\varphi (X_0, X_1)$ became
the well known Calder\'on spaces $X_0^{\theta } X_1^{1-\theta }$ (see \cite{Ca64}). Another important situation, investigated 
by Calder\'on and independently by Lozanovski{\u \i}, appears when $X_1\equiv L^{\infty }$. In particular,  
$X^{1/p}(L^{\infty})^{1-1/p} = X^{(p)}$ for $1 < p < \infty $ is known as the {\it $p$-convexification} of $X$ (see \cite{LT79}).

The properties of $\varphi (X_0, X_1)$ were studied by Lozanovski{\v \i} in \cite{Lo73, Lo78a} and \cite{Lo78b} (see also \cite{Ma89}), 
where among other facts it is proved the Lozanovski{\v \i} duality theorem: for any Banach ideal spaces $X_0, X_1$ 
with ${\rm supp} X_0 = {\rm supp} X_1$ and $\varphi \in {\mathcal U}$ we have
\begin{equation} \label{Lduality}
\varphi (X_0, X_1)^{\prime} = \hat{\varphi}(X_0^{\prime}, X_1^{\prime}),
\end{equation}
where the conjugate function $\hat {\varphi}$ is defined by
$$
\hat {\varphi} (s, t) := \inf \Big\{\frac {a s + b t} {\varphi (a, b)} ; \, a, b > 0\Big\}, s, t \geq 0.
$$
There hold $\hat {\varphi} \in {\mathcal U}$ and $\, \hat {{\hat \varphi}} = \varphi$ (see \cite[Lemma 2]{Lo78b}, \cite[Lemma 2]{Ma85} and 
\cite[Lemma 15.8]{Ma89}). It is easy to see that the Calder\'on-Lozanovski{\v \i} construction $\varphi (\cdot)$ is homogeneous 
with respect to an arbitrary weight $w$, that is, the equality
\begin{equation} \label{homogeneous}
\varphi (X_0(w), X_1(w)) = \varphi(X_0, X_1)(w),
\end{equation}
holds for arbitrary Banach ideal spaces $X_0, X_1$ and arbitrary weight $w$.

More information, especially on interpolation property, can be found in \cite{Be81, BK91, KLM13, KPS82, KMP93, Lo78b, Ma85, Ma89, 
Ni85, Ov76, Ov84, Sh81}.
\vspace{3mm}

We shall now identify the Calder\'on-Lozanovski{\v \i} construction for abstract Ces\`aro spaces.

%%%%%%%%%%%%%%%%%Theorem 2
\begin{theorem}\label{thm:inclusions}  
For any Banach ideal spaces $X_0, X_1$ and $\varphi \in {\mathcal U}$ the following embeddings hold
\begin{equation} \label{embeddings}
\varphi (CX_0,CX_1) \overset{1}{\hookrightarrow }C [\varphi(X_0,X_1)] ~~{\it and} ~~ \varphi (\widetilde{X_0},\widetilde{X_1})\overset{1}{\hookrightarrow }  [\varphi (X_0,X_1)]^{\Large \sim }.
\end{equation}
\end{theorem}
%%%%%%%%%%%%%%%%%%%%

%%%%  Proof
\proof
Suppose $X_0, X_1$ are Banach function spaces on $I = [0, \infty)$ and let $f \in \varphi(C{X_0},C{X_1})$ with 
$\|f \|_{\varphi(C{X_0},C{X_1})}<\lambda$. 
Then $ |f|\leq \lambda \varphi (|f_0|, |f_1|)$ a.e. for some $f_i\in CX_i$ with $\| f_i\|_{C{X_i}}\leq 1$ for $i=0,1$. 

Using now the Jensen inequality for concave function we obtain
\begin{eqnarray*}
C|f|(x) 
&=& 
\frac{1}{x}\int_0^x |f(t)| \, dt \leq \frac{\lambda}{x}\int_0^x \varphi (|f_0(t)|, |f_1(t)|) \, dt \\
&\leq& 
\lambda \varphi \Big(\frac{1}{x}\int_0^x |f_0(t)| \, dt, \frac{1}{x}\int_0^x |f_1(t)| \, dt \Big) =\lambda \varphi (C|f_0|(x), C|f_1|(x)).
\end{eqnarray*}
Since $\|C|f_i| \|_{{X_i}} = \| f_i\|_{C{X_i}}\leq 1$, $i=0,1$, it follows that $f \in C [\varphi({X_0},{X_1})]$ with 
$\|f\|_{C [\varphi({X_0},{X_1})]} = \|C|f| \|_{\varphi({X_0},{X_1})}\leq \lambda$ and the first embedding in (\ref{embeddings}) is proved.

To prove the second embedding in (\ref{embeddings}) let $f \in \varphi(\widetilde{X_0},\widetilde{X_1})$ with $\|f\|_{\varphi(\widetilde{X_0},\widetilde{X_1})}<\lambda$. This means that $ |f|\leq \lambda \varphi (|f_0|, |f_1|)$ a.e. for some $f_i\in \widetilde{X_i}$ with 
$\|f_i\|_{\widetilde{X_i}}\leq 1$,  $i=0,1$. Therefore,  for all $x \in I$,
\begin{eqnarray*}
\widetilde{f}(x) 
&=& 
\esssup_{t \in I, \, t \geq x} |f(t)| \leq \lambda \esssup_{t \in I, \, t \geq x} \varphi (|f_0(t)|, |f_1(t)|) \\
&\leq& 
\lambda \varphi \big (\esssup_{t \in I, \, t \geq x} |f_0(t)|, \esssup_{t \in I, \, t \geq x} |f_1(t)| \big) =\lambda \varphi (\widetilde{f_0}(x), \widetilde{f_1}(x)). 
\end{eqnarray*}
Of course, $\widetilde{f_i} \in X_i$ and $\|\widetilde{f_i}\|_{X_i}\leq 1$ for $i = 0, 1$, which means that $\widetilde{f}\in\varphi (X_0,X_1)$ 
or $f \in \widetilde{\varphi (X_0,X_1)}$ with $\|f\|_{\widetilde{\varphi (X_0,X_1)}}\leq \|f\|_{\varphi(\widetilde{X_0},\widetilde{X_1})}$. 
This proves the second embedding in (\ref{embeddings}) for the case of $I = [0, \infty)$. 

The remaining cases of $I = [0, 1]$ or $I = \mathbb N$ require only the evident modifications and therefore will be omitted.  
\endproof

Of course, the natural question is if there are equalities in (\ref{embeddings}) and, in fact, it is the case when we assume something more on spaces $X_0$ and $X_1$.

%%%%%%%%%%%%%%%%%%%%%% Theorem 3
\begin{theorem}\label{thm:cesCL}
Let $X_0, X_1$ be Banach ideal spaces with the Fatou property and such that the Ces\`aro operator $C$ is bounded on $X_0$ and $X_1$. Suppose that one of the following conditions hold:
\begin{itemize}
\item[(i)] $I = [0,\infty)$ and the dilation operator $\sigma_\tau$ is bounded in $X_0$ and $X_1$ for some $0 < \tau < 1$, 
\item[(ii)] $I = [0, 1]$ and $X_0, X_1$ are symmetric spaces with the Fatou property such that both operators 
$C, C^*: X_i \rightarrow X_i$ are bounded for $i =0, 1$, 
\item[(iii)] $I = \mathbb N$ and the dilation operator $\sigma_{3}$ is bounded on dual spaces $X_0^{\prime}$ and 
$X_1^{\prime}$.
\end{itemize}
Then 
\begin{equation} \label{3.4}
\varphi (CX_0, CX_1) = C [\varphi(X_0, X_1)].
\end{equation}
\end{theorem}
%%%%%%%%%%%%%%%

%%%%%%%%%
\proof
In view of Theorem \ref{thm:inclusions} we need to prove only the remaining inclusion. It appears however, that both inclusions 
in (\ref{embeddings}) are complemented to each other by duality. Thus we have for particular cases:

(i) Let $I = [0,\infty)$. Using twice the Lozanovski{\u \i} duality theorem (\ref{Lduality}), Theorem A(i) on duality of Ces\`aro spaces and the 
second imbedding from Theorem \ref{thm:inclusions} we obtain
\begin{eqnarray*}
[\varphi (CX_0, CX_1)]^{\prime}
&=&
\hat{\varphi}([CX_0]^{\prime}, [CX_1]^{\prime}) = \hat{\varphi}(\widetilde{X_0^{\prime}}, \widetilde{X_1^{\prime}}) 
\hookrightarrow [\hat{\varphi}(X_0^{\prime}, X_1^{\prime})]^{\thicksim}  \\
&=&
[{\varphi}(X_0, X_1)^{\prime}]^{\thicksim} = [C\varphi(X_0,X_1)]^{\prime}. 
\end{eqnarray*}
In the last equality, in order to use Theorem A(i), we notice that if $\sigma_{\tau}$ is bounded on $X_0$ and $X_1$ for some 
$0 < \tau < 1$, then it is also bounded on $\varphi(X_0,X_1)$ (one can also use here a more general result that $\varphi(X_0,X_1)$ is an interpolation space between $X_0$ and $X_1$ for positive operators -- see \cite[Theorem 1]{Be81}, \cite[Theorem 3.1]{Sh81}, \cite[Theorem 1]{Ma85} and \cite[Theorem 15.13]{Ma89}).

Finally, by the Fatou property of both spaces, we obtain that $CX_0, CX_1, \varphi (X_0, X_1)$ and $\varphi (CX_0, CX_1)$ have 
the Fatou property (cf.  \cite[Theorem 1(d)]{LM15a} and \cite[Corollary 3, p. 185]{Ma89}), and so
$$
C[\varphi(X_0, X_1)] \equiv [C \big(\varphi(X_0, X_1) \big)]^{\prime \prime} \hookrightarrow 
[\varphi (CX_0,CX_1)]^{\prime \prime} \equiv \varphi (CX_0,CX_1),
$$
which finishes the proof in this case. 

(ii) Let $I = [0, 1]$. Similarly as before, by the Lozanovski{\u \i} duality result (\ref{Lduality}) used twice, Theorem A(ii) on duality of 
Ces\`aro spaces, property (\ref{homogeneous}) and the second imbedding from Theorem \ref{thm:inclusions} we have 
\begin{eqnarray*}
[\varphi(CX_0, CX_1)]^{\prime}
&=&
\hat{\varphi}([CX_0]^{\prime},[CX_1]^{\prime}) =
\hat{\varphi} \big(\widetilde{X_0^{\prime}(1/v)},\widetilde{X_1^{\prime}(1/v)} \big)  \hookrightarrow
 [\hat{\varphi}(X_0^{\prime}(1/v), X_1^{\prime}(1/v))]^{\thicksim} \\
&=&
 [\hat{\varphi}(X_0^{\prime}, X_1^{\prime})(1/v)]^{\thicksim} = \big[{\varphi}(X_0, X_1)^{\prime}(1/v) \big]^{\thicksim} =
[C \big(\varphi(X_0,X_1) \big)]^{\prime}. 
\end{eqnarray*}
where the weight $v$ is $v(x)=1-x, x \in [0, 1)$.  
Observe that assumptions of the Theorem A(ii) are satisfied for ${\varphi}(X_0, X_1)$ thanks to interpolation property of the 
Calder\'on-Lozanovski{\v \i} construction for positive operators. Once again, by the Fatou property of both spaces, we have  
$$
C\big(\varphi(X_0, X_1)\big) \equiv [C \big(\varphi(X_0, X_1) \big)]^{\prime \prime} \hookrightarrow \varphi (CX_0,CX_1)^{\prime \prime} 
\equiv \varphi (CX_0,CX_1).
$$
Proof in the sequence case is similar to the proof of case (i).
\endproof

Immediately from Theorems \ref{thm:inclusions} and \ref{thm:cesCL} by taking $\varphi(s, t) = \min(s, t)$ and 
$\varphi (s, t) = \max (s, t) \approx s + t$ we obtain results for the intersection and the sum of Ces\`aro spaces. On the other hand,
taking $\varphi(s, t) = s^{1 - \theta} t^{\theta}, 0 < \theta < 1$ and $X_i = L^{p_i}$ or $X_i = l^{p_i}, i = 0, 1$ we obtain another corollary.

%%%%%%%%%%%% Corollary 1
\begin{corollary} \label{Co1}
Under the assumptions of Theorem 3 for Banach ideal spaces $X_0, X_1$ we have
\begin{equation} \label{Ces+}
C(X_0 \cap X_1) = C(X_0) \cap C(X_1) ~~ {\rm and} ~~ C(X_0 + X_1) = C(X_0) + C(X_1).
\end{equation}
\end{corollary}
Note that one can even prove that $C(X_0 \cap X_1) \equiv C(X_0) \cap C(X_1)$ without additional assumptions on the spaces, but
the second equality in (\ref{Ces+}) is not true in general. If we consider spaces on $[0, \infty)$, then 
$C(L^1) + C(L^{\infty}) = C(L^{\infty})\subsetneq C(L^1 + L^{\infty})$.

%%%%%%%%%%%%% Example 2
\begin{example} \label{Ex2} 
Let $f(x) = x^{-1} \ln^{-3} (1/x) \, \chi_{(0, 1/e)}(x)$ for $x > 0$. Then $\| f \|_{L^1} = 1/2$ and 
$f \in C(L^1 + L^{\infty}) \setminus C(L^{\infty}) $. In fact, since
$Cf(x) = \frac{1}{2 x \ln^2x} \, \chi_{(0, 1/e]}(x) + \frac{1}{2 x} \, \chi_{(1/e, \infty)}(x)$ it follows that
$$
\| f \|_{C(L^{\infty}) } = \sup_{x > 0} Cf(x) = \lim_{x \rightarrow 0^+}  \frac{1}{2 x \ln^2x}  = \infty,
$$
and
\begin{eqnarray*}
\| f \|_{C(L^1 + L^{\infty}) } =
\int_0^1 \big( C|f|)^*(x) \, dx = \int_0^{1/e}  \frac{1}{2 x \ln^2x} \, dx + \int_{1/e}^1 \frac{1}{2 x} \, dx= 1.
\end{eqnarray*}
It seems to be of interest to investigate structure of the space $C(L^1 + L^{\infty})$ with its norm 
$\| f \|_{C(L^1 + L^{\infty})} = \int_0^1 \big( C|f|)^*(x) \, dx$.
\end{example}

%%%%%%%%%%%% Corollary 2
\begin{corollary} \label{Co2}
Let $1< p_0, p_1<\infty$ and $0<\theta<1$ be such that $\frac{1}{p}=\frac{1-\theta}{p_0}+\frac{\theta}{p_1}$. Then 
\begin{equation} \label{Cestheta}
\big( Ces_{p_0} \big)^{1-\theta} \big(Ces_{p_1} \big)^{\theta} = Ces_{p} ~~ {\rm and} ~~ ( ces_{p_0})^{1-\theta} (ces_{p_1})^{\theta} = ces_{p}.
\end{equation}
\end{corollary}

Note that by the duality Theorem A and Theorem \ref{thm:cesCL} we have also equality 
$\varphi(\widetilde{X_0}, \widetilde{X_1}) = [\varphi(X_0, X_1)]^{\Large \sim}$. It is however instructive to see that assumptions on  $X_0$
 and $X_1$ may be weakened if we do not use the duality argument in the proof. 

%%%%%%%%%%%%%%%%%%%%%%%% Theorem 4
\begin{theorem}\label{thm:tyldaspace}
Let $X_0$ and $X_1$ be two Banach ideal spaces with the Fatou property. If either $X_0, X_1$ are symmetric spaces or
$C, C^*$ are bounded on both $X_0$ and $X_1$, then 
\begin{equation} \label{tilda}
\varphi(\widetilde{X_0}, \widetilde{X_1}) = [\varphi(X_0, X_1)]^{\Large \sim}.
\end{equation}
\end{theorem}
%%%%%%%%%%%%%%

%%%%%%%%
\proof
Of course, we need to prove only the inclusion ``$\hookrightarrow$" and the Fatou property is not necessary here. We restrict ourselfs to the case of function spaces on $I=[0,\infty)$, since the remaining cases are analogous. Let 
$f \in \widetilde{\varphi(X_0, X_1)}$ with $\| f \|_{\widetilde{\varphi(X_0, X_1)}}<\lambda$. Then 
$\widetilde{f}\in \varphi(X_0,X_1)$ and so there are $f_i\in X_i$ with $\| f_i\|_{X_i}\leq 1$ for $i=0,1$ such that 
$$
\widetilde{f}(x) \leq \lambda \varphi (|f_0(x)|, |f_1(x)|) ~~{\rm a.e. ~on} ~I.
$$
Now, if we assume that both spaces $X_0$ and $X_1$ are symmetric, then the argument is as follows. Firstly, we have
$$
\widetilde{f}(x) = \big (\widetilde{f} \, \big)^*(x) \leq \lambda \varphi (f_0^*(x/2), f_1^*(x/2)) ~ {\rm a.e. ~on} ~I,
$$
where the last inequality is a consequence of the estimation
\begin{equation} \label{inequality-star}
\varphi (|f_0|, |f_1|)^*(x) \leq \varphi (f_0^*(x/2), f_1^*(x/2)) ~~{\rm for} ~x \in I.
\end{equation}
Observe that the inequality (\ref{inequality-star}) follows from the fact that the conjugation operation on ${\mathcal U}$ is 
an  involution. In fact, by a standard inequality for rearrangements $(f+g)^*(x) \leq f^*(x/2) + g^*(x/2)$ we obtain
$$
\varphi (|f_0(x)|, |f_1(x)|) \leq \dfrac{a |f_0(x)| + b |f_1(x)|}{\hat{\varphi}(a, b)}~~{\rm for ~all} ~~a, b > 0
$$ 
and
$$
\varphi (|f_0|, |f_1|)^*(x) \leq \dfrac{a f_0^*(x/2) + b f_1^*(x/2)}{\hat{\varphi}(a, b)} ~~{\rm for ~all} ~~a, b > 0.
$$ 
Taking the infimum over all $a, b > 0$ we get
$$
\varphi (|f_0|, |f_1|)^*(x) \leq  \hat {{\hat \varphi}}  (f_0^*(x/2), f_1^*(x/2)) = \varphi (f_0^*(x/2), f_1^*(x/2)).
$$ 
Secondly, putting $g_i(x) = f_i^*(x/2)$ for $i = 0, 1$, by symmetry of $X_i$, we obtain $g_i \in X_i$ and 
$$
\| g_i\|_{X_i} \leq \|\sigma_2\|_{X_i\rightarrow X_i}\| f_i\|_{X_i} \leq \|\sigma_2\|_{X_i\rightarrow X_i}, i = 0, 1.
$$ 
Moreover, $\widetilde{g_i} = g_i$ and so $\| g_i\|_{\widetilde{{X_i}}} \leq \|\sigma_2\|_{X_i\rightarrow X_i}$
for $i = 0, 1$. Thus  $f \in \varphi(\widetilde{X_0}, \widetilde{X_1})$ with $\| f\|_{\varphi(\widetilde{X_0}, \widetilde{X_1})}\leq A \|f\|_{\widetilde{\varphi(X_0,X_1)}}$, where 
$A = \max (\| \sigma_2\|_{X_0\rightarrow X_0}, \| \sigma_2\|_{X_1\rightarrow X_1})$.
\vspace{3mm}

If we assume that both $C$ and $C^*$ are bounded on $X_0$ and $X_1$, then we use the following estimation
\begin{equation} \label{3.9}
CC^*[ \varphi(|f_0|, |f_1|)] \leq \varphi[CC^*(|f_0|), CC^*(|f_1|)] ~ {\rm a.e. ~ on } ~ I.
\end{equation}  
In fact, we have for a.e. $x \in I$ and all $a, b > 0$ 
$$
\varphi(|f_0(x)|, |f_1(x)|) \leq \dfrac{a |f_0(x)| + b |f_1(x)|}{\hat{\varphi}(a, b)}.
$$
Then, by linearity and monotonicity of $C$ and $C^*$ we get for $x \in I$
$$
CC^*[\varphi(|f_0|, |f_1|)](x) \leq \frac{a \, CC^*|f_0(x)| + b \, CC^*|f_1(x)|}{\hat{\varphi}(a, b)}
$$
for all $a,b>0$. Using again the involution property of the conjugation we obtain
$$
CC^*[\varphi(|f_0|, |f_1|)] \leq \hat{\hat{\varphi}}(CC^*|f_0|, CC^*|f_1|) = \varphi(CC^*|f_0|, CC^*|f_1|).
$$
Consequently, 
\begin{equation} \label{3.10}
\begin{split}
|f| &
\leq 
\widetilde{f}\leq C\widetilde{f}\leq C\widetilde{f}+C^*\widetilde{f}= CC^*\widetilde{f} \\
&\leq 
CC^*[\varphi(|f_0|, |f_1|)] \leq \varphi[CC^*(|f_0|), CC^*(|f_1|)]. \\
\end{split}
\end{equation}
However, by our assumption $CC^*(|f_i|) \in X_i$ and $\widetilde{CC^*|f_i|} = CC^*|f_i|$, which means that 
$CC^*|f_i| \in \widetilde{X_i}$ for $i = 0, 1$. Therefore, 
$$
f \in \varphi(\widetilde{X_0}, \widetilde{X_1}) ~ ~{\rm and} ~~ \|f \|_{\varphi(\widetilde{X_0}, \widetilde{X_1})}\leq B \| f \|_{\widetilde{\varphi(X_0, X_1)}},$$
where $B = \max (\|C\|_{X_0\rightarrow X_0}\|C^*\|_{X_0\rightarrow X_0}, \|C\|_{X_1\rightarrow X_1}\|C^*\|_{X_1\rightarrow X_1})$. 
\endproof

%%%%%%%%%%%% Corollary 3
\begin{corollary}  \label{Co3}
Let $1< p_0, p_1<\infty$ and $0<\theta<1$ be such that $\frac{1}{p}=\frac{1-\theta}{p_0}+\frac{\theta}{p_1}$. Then 
\begin{equation} \label{Cestheta}
\big( \widetilde{L^{p_0}} \big)^{1-\theta} \big( \widetilde{L^{p_1}} \big)^{\theta} =  \widetilde{L^p} ~~ {\rm and} ~~ ( \widetilde{l^{p_0}} )^{1-\theta} (\widetilde{l^{p_1}} )^{\theta} = \widetilde{l^{p}}.
\end{equation}
\end{corollary}

Applying Theorem \ref{thm:tyldaspace} we can get some ``one-sided" result similar to (\ref{3.4}).

%%%%%%%%%%%% Corollary 4
\begin{corollary} \label{Co4} (a) Let $I = [0, \infty)$. If $X$ is a symmetric space on $I$ with the Fatou property such that $C$ is bounded 
on $X$ and on $\varphi (L^1, X)$, then
\begin{equation} \label{3.12}
\varphi (L^1, CX) = C[\varphi (L^1, X)].
\end{equation}
In particular,
\begin{equation} \label{3.13}
(L^1)^{1-\theta} (Ces_{p})^{\theta} = Ces_{q} 
\end{equation}
for $1 < p \leq \infty, \frac{1}{q} = 1- \theta + \frac{\theta}{p}$ and any $0 < \theta < 1$.

(b) Let $I = [0, 1]$. Then
\begin{equation} \label{3.14}
\varphi (L^1, Ces_{\infty}) = C[\varphi (L^1, L^{\infty})].
\end{equation}
\end{corollary}

%%%%%%%%
\proof
(a) Using twice the Lozanovski{\v \i} duality result, twice Theorem A and Theorem \ref{thm:tyldaspace} 
we obtain
\begin{eqnarray*}
\varphi(L^1, CX)^{\prime} 
&=&
\hat{\varphi}[(L^1)^{\prime}, (CX)^{\prime}] = \hat{\varphi}(L^{\infty}, \widetilde{X^{\prime}}) 
= \hat{\varphi}( \widetilde{L^{\infty}}, \widetilde{X^{\prime}}) \\
&=&
\widetilde{\hat{\varphi}( L^{\infty}, X^{\prime}) } = \widetilde{\hat{\varphi}( (L^1)^{\prime}, X^{\prime})} = 
\widetilde{\varphi ( L^1, X)^{\prime}} = [C \varphi (L^1, X)]^{\prime}.
\end{eqnarray*}
Then, by the Fatou property of $X$, all $CX, \varphi (L^1, CX), \varphi (L^1, X)$ and $C[\varphi (L^1, X)]$ also have 
 the Fatou property (cf. \cite[Theorem 1(d)]{LM15a} and \cite[Corollary 3, p. 185]{Ma89}), and so
$$
\varphi(L^1, CX) \equiv \varphi(L^1, CX)^{\prime \prime} = [C \varphi (L^1, X)]^{\prime \prime} \equiv C \varphi (L^1, X).
$$
Equality (\ref{3.13}) follows from (\ref{3.12}) and the identification $(L^1)^{1-\theta} (L^p)^{\theta} = L^q$.

(b) Similarly as in (a), with Theorem A replaced by the Luxemburg-Zaanen duality result \cite{LZ66} (see also 
\cite[Theorem 7]{LM15a}) $(Ces_{\infty})^{\prime}  \equiv \widetilde{L^1}$, we obtain
\begin{eqnarray*}
\varphi(L^1, Ces_{\infty})^{\prime} 
&=&
\hat{\varphi}[(L^1)^{\prime}, (Ces_{\infty})^{\prime}] = \hat{\varphi}(L^{\infty}, \widetilde{L^1}) 
= \hat{\varphi}( \widetilde{L^{\infty}}, \widetilde{L^1}) \\
&=&
\widetilde{\hat{\varphi}( L^{\infty}, L^1) } = \widetilde{\hat{\varphi}( (L^1)^{\prime}, (L^{\infty})^{\prime})} = 
\widetilde{\varphi ( L^1, L^{\infty})^{\prime}} = [C \varphi (L^1, L^{\infty})]^{\prime},
\end{eqnarray*}
and by the Fatou property $\varphi(L^1, Ces_{\infty}) = C \varphi (L^1, L^{\infty})$.
\endproof

%%%%%%%%%%%%%%%%% Example 3
\begin{example} \label{Ex3}
In the case $I = [0,1]$ one cannot expect a general result like in Corollary \ref{Co4}(a) even for $X = L^2$. In fact,
for the weight $v(x) = 1 - x$ we have 
$$
[(L^1)^{1/2}(Ces_2)^{1/2}]^{\prime} = (L^{\infty})^{1/2}(\widetilde{L^2(1/v)})^{1/2} = (\widetilde{L^2(1/v)})^{(2)}=\widetilde{L^4(1/\sqrt{v})},
$$
where $X^{(2)}$ is $2$-convexification of $X$. On the other hand, 
$$
(C[(L^1)^{1/2}(L^2)^{1/2}])^{\prime} = [C(L^{4/3})]^{\prime} = (Ces_{4/3})^{\prime} = \widetilde{L^4(1/v)}.
$$ 
Therefore, $(L^1)^{1/2}(Ces_2)^{1/2} = [(L^1)^{1/2}(Ces_2)^{1/2}]^{\prime \prime} = [\widetilde{L^4(1/\sqrt{v})}]^{\prime}$ and
$$
C[(L^1)^{1/2}(L^2)^{1/2}] = (C[(L^1)^{1/2}(L^2)^{1/2}])^{\prime \prime} [ \widetilde{L^4(1/v)}]^{\prime} = Ces_{4/3} \subsetneq  
[\widetilde{L^4(1/\sqrt{v})}]^{\prime}
$$
since, of course, $\widetilde{L^4(1/v)} \subsetneq \widetilde{L^4(1/\sqrt{v})}$. 
\end{example}

Since we deal with Banach ideal spaces, the previous results one can apply to the complex method of interpolation. In order to 
present it we need the following simple lemma. 

%%%%%%%%%%%%%%%%%%%% Lemma 1
\begin{lemma} \label{L1}  \item[(a)] If a Banach ideal space $X$ satisfies $X \in (OC)$, then $CX \in (OC)$. 
\item[(b)] If a Banach sequence space $X\subset c_0$ satisfies $X \in (OC)$, then $\widetilde{X} \in (OC)$. 
\end{lemma}
%%%%%%%%%%%%%%%%%%
\proof 
(a) Let $f\in CX$, where $X$ is a Banach ideal space and let $(A_n)$ be a sequence of measurable sets such that 
$A_{n+1}\subset A_n$, $n=1,2,3,...$, and $m(\bigcap A_n) = 0$. Then, by the Lebesgue domination theorem, for each $x > 0$ 
$$
C|f|(x) = \frac{1}{x} \int_0^x |f(t)|\chi_{A_n}(t) \, dt \rightarrow 0 \ {\rm as }\ n\rightarrow \infty.
$$
Thus, $C|f|\geq C|f|\chi_{A_n}\rightarrow 0$ a.e. and by order continuity of $X$ we obtain $\| f \chi_{A_n}\|_{CX} 
= \|C|f|\chi_{A_n}\|_X\rightarrow 0$, as required. Simple modification proves the case of sequence spaces.

(b) Let $x\in \widetilde{X}$. It is enough to check whether $\|x\chi_{[n,\infty)}\|_{\widetilde{X}}\rightarrow 0$. We have 
$\tilde x\geq \widetilde{x\chi_{[n,\infty)}}$. On the other hand, since  $X\subset c_0$, one gets $ (\widetilde{x\chi_{[n,\infty)}})(i)\leq \widetilde{x}_n\rightarrow 0$ for each $i=1,2,3,\dots$, so that order continuity of $X$ gives the claim. 
\endproof

Notice that the reverse implication in the above Lemma \ref{L1}(a) does not hold (cf. Example \ref{Ex1}):  for 
$X = L^2[0, 1/4] \oplus L^{\infty}[1/4, 1/2] \oplus L^2[1/2, 1]$ we have $CX[0, 1] = CL^2[0, 1] = Ces_2[0, 1] \in (OC)$ but $X\not \in (OC)$. 
It is also worth to emphasize here once more conclusion of Lemma \ref{L1}(b). Namely, in contrast to the Tandori function spaces, which are never order continuous (cf. \cite{LM15a}), the construction in sequence case behaves quite well.
\vspace{2mm}

To simplify an exposition of the next theorem we assume that Banach ideal spaces contain all characteristic functions of subsets of finite measure of 
underlying measure space. This assumption is equivalent to the imbedding (cf. \cite[Lemma 4.1, p. 90]{KPS82}):
\begin{equation} \label{Eq3.15}
L^1 \cap L^{\infty} \hookrightarrow X.
\end{equation}
Condition (\ref{Eq3.15}) is a little stronger than existence of a weak unity. Moreover, if $X$ is a Banach ideal space and (\ref{Eq3.15}) 
holds, then also $\widetilde{X}$ is a Banach ideal space with property (\ref{Eq3.15}). If additionally $C$ is bounded on $X$, then also 
$CX$ satisfies property (\ref{Eq3.15}).
\vspace{2mm}

%%%%%%%%%%%%%%%%%%%%%%%%%%%%%%%%%%
The {\it Calder\'on (lower) complex method of interpolation} is defined only for a couple of Banach spaces $(X_0, X_1)$ over 
the complex field therefore we must, in fact, apply it to the couple $(X_0({\mathbb C}), X_1({\mathbb C}))$, where $X_k({\mathbb C})$ 
denotes the {\it complexification} of $X_k$ (namely the space of all complex-valued measurable functions $f: I \rightarrow \mathbb C$ 
such that $|f| \in X_k$ with the norm $\| f\|_{X_k({\mathbb C})} = \| |f| \|_{X_k}), k = 0, 1$. Let $[X_0, X_1]_{\theta}$ denote the subspace 
of real-valued functions in Calder\'on's interpolation spaces $[X_0({\mathbb C}), X_1({\mathbb C})]_{\theta}$ for $0 < \theta < 1$. For 
formal definition and properties of the Calder\'on (lower) method of complex interpolation we refer to original Calder\'on's paper \cite{Ca64} 
and books \cite{BL76, BK91, KPS82}.

All proofs given for Banach ideal spaces of real-valued functions are true also for complexified Banach ideal spaces of measurable 
functions on $I$ (cf. \cite{Ca64, CN03, Cw10}). For example, if ${\rm supp}X_0 = {\rm supp}X_1 = I$, then 
${\rm supp}X_0 \cap X_1 = {\rm supp}[X_0, X_1]_{\theta} = I$.

%%%%%%%%%%%%%%%%%%%%%%%%%
\begin{theorem}\label{complex}
Let $0 < \theta < 1$. Assume that $X, X_0, X_1$ are Banach ideal spaces on $I$ with the Fatou property and the property (\ref{Eq3.15}), 
and such that the Ces\`aro operator $C$ is bounded on all of them. 
\begin{itemize}
\item[(a)] If $I = [0, \infty)$, the dilation operator $\sigma_a$ for some $0<a<1$ is bounded on $X_0$ and $X_1$ and at least one of the spaces 
$X_0$ or $X_1$ is order continuous, then 
$$
[CX_0, CX_1]_{\theta} = C([X_0, X_1]_{\theta}).
$$
\item[(b)]  If $I = [0, 1]$, $X_0, X_1$ are symmetric spaces such that $C^*$ is bounded on both of them and at least one of the spaces 
$X_0$ or $X_1$ is order continuous, then 
$$
[CX_0, CX_1]_{\theta} = C([X_0, X_1]_{\theta}).
$$
\item[(c)] Let $X_0,X_1$ be  Banach sequence spaces such that the dilation operator $\sigma_{3}$ is bounded on dual spaces 
$X_0^{\prime}, X_1^{\prime}$ and at least one of the spaces $X_0$ or $X_1$ is order continuous, then 
$$
[CX_0,CX_1]_{\theta} = C([X_0,X_1]_{\theta}).
$$
\item[(d)] If $X$ is a symmetric space on $I=[0,\infty)$, then
$$
[L^1, CX]_{\theta} = C([L^1, X]_{\theta}).
$$
\item[(e)] For $I = [0, 1]$ we have
$$
[L^1, Ces_{\infty}]_{\theta} = C([L^1, L^{\infty}]_{\theta}).
$$
\item[(f)] Let  $I = [0, \infty)$ or  $I = [0, 1]$ and suppose that at least one of the spaces $X_0,X_1$ is order continuous. If either $X_0$ and 
$X_1$ are symmetric spaces or $C^*$ is bounded on $X_0$ and $X_1$, then 
$$
[\widetilde{X_0},\widetilde{X_1}]_{\theta} = ([X_0,X_1]_{\theta})^{\sim}.
$$
\end{itemize}
\end{theorem}
%%%%%%%%%%%%%%%%%%%%%%
\proof The main tool in the proof of all points will be Theorem \ref{thm:cesCL} and Shestakov's representation of the complex method of interpolation for Banach ideal spaces $X_0, X_1$ (cf. \cite[Theorem 1]{Sh74}; see also \cite[Theorem 9]{RT10}), i.e.,  
\begin{equation} \label{Shestakov}
[X_0, X_1]_{\theta} \equiv \overline{X_0\cap X_1}^{X_0^{1-\theta}X_1^{\theta}}.
\end{equation}
In fact, the proofs of all our points relay on this theorem. It is enough to notice that in all cases from (a) to (c) we have 
$$
\overline{X_0\cap X_1}^{X_0^{1-\theta}X_1^{\theta}} = {X_0^{1-\theta}X_1^{\theta}} ~~ {\rm and} ~~  \overline{CX_0\cap CX_1}^{CX_0^{1-\theta}CX_1^{\theta}} = ({CX_0)^{1-\theta} (CX_1)^{\theta}}
$$
just because $X_0^{1-\theta}X_1^{\theta}$ and $(CX_0)^{1-\theta} (CX_1)^{\theta}$ are order continuous under our assumptions. In fact, 
$X_0^{1-\theta} X_1^{\theta}$ for $0 < \theta < 1$ is order continuous when at least one of $X_0$ or $X_1$ is order continuous (see, for 
example, \cite[Lemma 20, p. 428]{Lo69}, \cite[Proposition 4]{Re88}, \cite[Theorem 15.10]{Ma89} and \cite{KL10}, where $(OC)$ property 
of $\varphi(X_0, X_1)$ was investigated), so that simple functions are dense therein. Using all of these representations we get 
\begin{equation*}
[CX_0,CX_1]_{\theta} = C([X_0,X_1]_{\theta}) = C(X_0^{1-\theta}X_1^{\theta}) = (CX_0)^{1-\theta} (CX_1)^{\theta}
\end{equation*}
with a corresponding modifications in points (d) and (e), where instead of Theorem \ref{thm:cesCL} we use Corollary \ref{Co4}.

To prove (f) we need a little more delicate argument. First of all recall that $\widetilde{X}$ are never order continuous in a function case 
and even  worse $\widetilde{X}_a = \{0\}$ (see \cite[Theorem 1(e)]{LM15a}). So that at the first look it seems to be hopeless to apply an argument like above here. Fortunately, order continuity means only that one can approximate any function in a norm by each majorized sequence tending to it almost everywhere, but we need to approximate a given function at least by one sequence, so that the lack of order continuous elements will not be an obstacle. Of course,
\begin{equation*}
\overline{\widetilde{X_0}\cap\widetilde{X_1}}^{\widetilde{X_0}^{1-\theta}\widetilde{X_1}^{\theta}} \hookrightarrow \widetilde{X_0}^{1-\theta}\widetilde{X_1}^{\theta}.
\end{equation*}
We will show that in the above imbedding we have equality. Firstly, consider the case when $I = [0,\infty)$. Let 
$f\in \widetilde{X_0}^{1-\theta}\widetilde{X_1}^{\theta}=\widetilde{X_0^{1-\theta}X_1^{\theta}}$ be such that $\widetilde{f}(x) > 0$ for 
each $x \geq 0$ (always $\widetilde{f}(x) \geq 0$ and we do proof in the worst case when $\widetilde{f}(x) > 0$). By definition $\tilde f\in X_0^{1-\theta}X_1^{\theta}$. Since $\tilde f$ is nonincreasing, for each interval $[\frac{1}{n},n]$ we can find a simple function $g_n$ with support in  $[\frac{1}{n},n]$ such that, $|g_n|\leq|f|$ and $\|(f-g_n)\chi_{[\frac{1}{n},n]}\|_{\infty} \leq \widetilde{f}(n)$.  
Clearly, $g_n\in \widetilde{X_0}\cap\widetilde{X_1}$. We have 
$$
| f - g_n |\leq |f|\chi_{[0,\frac{1}{n}]}+\widetilde{f}(n)\chi_{[\frac{1}{n},n]}+|f|\chi_{[n,\infty)}
$$
and consequently
$$
\widetilde{|(f-g_n)|}\leq \widetilde{|f|\chi_{[0,\frac{1}{n}]}}+\widetilde{f}(n)\chi_{[0,n]}+\widetilde{|f|\chi_{[n,\infty)}}.
$$
Now, since $\widetilde{f}(n) \rightarrow 0$ with $n\rightarrow \infty$, we see that $ \widetilde{|f|\chi_{[0,\frac{1}{n}]}}\rightarrow 0$,  
$ \widetilde{f}(n) \chi_{[0,n]} \rightarrow 0$ and  $ \widetilde{|f|\chi_{[n,\infty)}}\rightarrow 0$ a.e. on $I$. Moreover, all these three 
sequences are dominated by $\widetilde{f}$, so that order continuity of $X_0^{1-\theta}X_1^{\theta}$ guarantees that 
$
\|\widetilde{f-g_n}\|_{X_0^{1-\theta}X_1^{\theta}}\rightarrow 0$ as $n \rightarrow \infty$, 
which proves the claim. In case when there is $a> 0$ such that $\widetilde{f}(x)=0$ for all $x > a$, we can proceed analogously, only
replacing intervals $[\frac{1}{n},n]$ by $[\frac{1}{n},n]\cap [0,a]$, $\widetilde{f}(n)$ by $\frac{1}{n}$ and a new majorant is then 
$\max\{\tilde f, \chi_{[0,a]}\}$. The same argument works as well for the case $I=[0,1]$.
\endproof

%%%%%%%%%%%%%%%%%%%%%%%%% Section 3
\section{Real method}

One of the most important interpolation methods is the {\it $K$-method} known also as the {\it real Lions-Peetre interpolation method}. For a Banach couple $\bar{X} = (X_0, X_1)$ the {\it Peetre K-functional} of an element $f \in X_0+X_1$ is defined for $t > 0$ by
$$
K(t, f; X_0, X_1) = \inf \{ \| f_0\|_{X_0} + t \| f_1\|_{X_1}: f = f_0 + f_1, f_0 \in X_0, f_1 \in X_1 \}.
$$
Let $G$ be a Banach ideal space on $(0, \infty)$ containing function $\min(1, t)$. Then the space of real interpolation 
or the $K$-method of interpolation
$$
(X_0, X_1)_G^K = \{f \in X_0 + X_1: K(t, f; X_0, X_1) \in G \}
$$
is a Banach space with the norm $\| f \|_{(X_0, X_1)_G^K} = \| K(t, f; X_0, X_1) \|_G$. This space is an intermediate space 
between $X_0$ and $X_1$, that is, $X_0 \cap X_1 \hookrightarrow (X_0, X_1)_G^K \hookrightarrow X_0 + X_1$. Moreover, $
(X_0, X_1)_G^K$ is an interpolation space between $X_0$ and $X_1$. 

The most common with several applications is $K$-method, where $G$ is given by 
$\| f\|_G = (\int_0^{\infty} (t^{-\theta} |f(t)|)^p \frac{dt}{t})^{1/p}$ for $0 < \theta < 1, 1 \leq p < \infty$ or 
$\| f\|_G = (\sup_{t > 0} t^{-\theta} |f(t)|)$ for $0 \leq \theta \leq 1$ (when $p= \infty$) and then 
$ (X_0, X_1)_{G}^K = (X_0, X_1)_{\theta, p}$.
More information about interpolation spaces, and, in particular, interpolation functors may be found in the books \cite{BS88, BL76, BK91} and \cite{KPS82}.

%%%%%%%%%%%%%%%%%%%%%%%% Theorem 6
\begin{theorem}\label{thm:Cesaro=real}
Let $X_0,X_1$ be two Banach function spaces on $I = [0, \infty)$. If $C$ and $C^*$ are bounded on $X_i$ for $i = 0, 1$ and $F$ is an interpolation functor 
with the homogenity property, that is, $F(X_0(w), X_1(w)) = F(X_0, X_1)(w)$ for any weight $w$ on $I$, then
\begin{equation} \label{Cesaro=real}
F(CX_0, CX_1) = CF(X_0, X_1).
\end{equation}
In particular, 
\begin{equation} \label{4.2}
(CX_0, CX_1)_G^K = C[(X_0, X_1)_G^K].
\end{equation}
\end{theorem}
%%%%%%%%%%%%%%

%%%%%%%%
\proof
Let $X(w_0)$ be a weighted Banach ideal space on $(0, \infty)$ with the weight $w_0(t) = \frac{1}{t}$ and such that $C, C^*$ are 
bounded on $X$. First of all notice that 
$$
K(t, f; L^1, L^1(1/s)) = K(t, |f|; L^1, L^1(1/s)) = t \, [C|f|(t) + C^*|f|(t)] 
$$
implies
\begin{equation} \label{CesaroK}
(L^1, L^1(1/s))_{X(w_0)}^K = CX.
\end{equation} 
In fact, we have
$$
\| f \|_{(L^1, L^1(1/s))_{X(w_0)}^K } = \| K(t, f; L^1, L^1(1/s)) \|_{X(w_0)}^K = \| C|f| + C^*|f| \|_X
$$
and 
\begin{eqnarray*}
\| f \|_{CX} 
&=& 
\| C|f| \|_X \leq \| C|f| + C^*|f|\|_X = \| C^*C|f| \|_X \\
&\leq& 
\| C^* \|_{X \rightarrow X} \, \| C|f| \|_X = \| C^* \|_{X \rightarrow X} \, \| f \|_{CX}.
\end{eqnarray*}

Notice now that the operator $Sf(t) = \int_0^{\infty} \min(1, \frac{t}{s}) f(s) \frac{ds}{s}$ is bounded on $X(w_0)$. In fact, 
\begin{eqnarray*}
\| Sf \|_{X(w_0)} 
&=&
\| \frac{1}{t} \int_0^t f(s) \frac{ds}{s} + \int_t^{\infty} \frac{f(s)}{s} \frac{ds}{s} \|_X \\
&=&
\| C(f w_0) + C^*(f w_0) \|_X \leq \big( \| C \|_{X \rightarrow X} +  \| C^* \|_{X \rightarrow X} \big) \| f w_0\|_X \\
&=&
 \big( \| C \|_{X \rightarrow X} +  \| C^* \|_{X \rightarrow X} \big) \| f \|_{X(w_0)}.
\end{eqnarray*}
Let's return to the proof of commutativity (\ref{Cesaro=real}). Since $C, C^*$ are bounded on $X_0$ and on $X_1$, then by (\ref{CesaroK})
we obtain $CX_i = (L^1, L^1(1/s))_{X_i(w_0)}^K, i = 0, 1$ and since operator $S$ is bounded on $X_0(w_0)$ and on $X_1(w_0)$, 
then by Brudny{\v \i}-Dmitriev-Ovchinnikov theorem (see \cite[Theorem 1]{DO79} and \cite[Theorem 4.3.1]{BK91}): for any interpolation 
functor $F$ we have 
$$
F \big( (L^1, L^1(1/s))_{X_0(w_0)}^K, (L^1, L^1(1/s))_{X_1(w_0)}^K \big) = (L^1, L^1(1/s))_{F(X_0(w_0), X_1(w_0))}^K.
$$ 
Now, by assumption that the functor $F$ has the homogeneity property we obtain
\begin{eqnarray*}
F(CX_0, CX_1) 
&=&
F \big( (L^1, L^1(1/s))_{X_0(w_0)}^K, (L^1, L^1(1/s))_{X_1(w_0)}^K \big) \\
&=&
(L^1, L^1(1/s))_{F(X_0(w_0), X_1(w_0))}^K \\
&=&
(L^1, L^1(1/s))_{F(X_0, X_1)(w_0)}^K = C[F(X_0, X_1)],
\end{eqnarray*}
where the last equality follows from (\ref{CesaroK}) and the fact that $F$ is an interpolation functor which implies boundedness 
of $C$ and $C^*$ on $F(X_0, X_1)$.

To prove the second statement, note that the $K$-method of interpolation is homogeneous. It follows from the equality 
$K(t, f; X_0(w), X_1(w)) = K(t, f w; X_0, X_1)$ (cf. \cite[Proposition 14]{LM15a}), since then
\begin{eqnarray*}
\| f \|_{(X_0(w), X_1(w))_G^K} 
&=& 
\| K(t, f; X_0(w), X_1(w)) \|_G = \| K(t, f w; X_0, X_1) \|_G \\
&=& 
\| f w \|_{(X_0, X_1)_G^K} = \| f  \|_{(X_0, X_1)_{G}^{K} (w)}. 
\end{eqnarray*}
\endproof

%%%%%%%%%%% Remark 2
\begin{remark}
Note that (\ref{4.2}) for the case of symmetric spaces $X_0, X_1$ on $[0, \infty)$ with the operator $C$ bounded on $X_i, i = 0, 1$  was already proved 
in \cite[Corollary 3.2]{AM13} and (\ref{CesaroK}) for the case $X = L^p$ with $1 < p < \infty$  was proved in \cite[Theorem 2.1(ii)]{AM13}. 
\end{remark}
%%%%%%%%%%%%%%% Remark 3
\begin{remark}
As it was mentioned in the previous section also the Calder\'on-Lozanovski{\v \i} construction has the homogeneity property (note that it is an interpolation functor for example when all spaces have the Fatou property). Another classical functors with the homogeneity property are functors of orbit $Orb_{\bar E}(a,\bar X,L)$ and coorbit $Corb_{\bar Y}(F,\bar X)$ (cf. \cite{BK91}). Using a similar argument as in the proof of Theorem \ref{complex} one could also prove the homogeneity of the complex method.
\end{remark}

As an example we consider interpolation of weighted Ces\'aro spaces $Ces_{p, \alpha} = C(L^p(x^{\alpha}))$ on $I = [0, \infty)$. We only 
need to observe that $C$ and $C^*$ are bounded on $L^p(x^{\alpha})$ if and only if $1 \leq p \leq \infty$ and $-1/p < \alpha < 1- 1/p$ 
(see \cite[p. 245]{HLP52} for sufficiency and \cite[pp. 38-40]{KMP07} for equivalence).

%%%%%%%%%%%% Corollary 5
\begin{corollary} \label{Co5}
Let $1 \leq p_i \leq\infty$ and $-1/p_i <\alpha_i <1-1/p_i$ for $i = 0, 1$, then 
\begin{equation} \label{Cesthetap}
\big ( Ces_{p_0, \alpha_0}, Ces_{p_1, \alpha_1} \big)_{\theta, p} = Ces_{p, \alpha},
\end{equation}
where $\frac{1}{p} = \frac{1-\theta}{p_0} + \frac{\theta}{p_1}$ and $\alpha = (1-\theta) \alpha_0 + \theta \alpha_1$
\end{corollary}

In particular, for $\alpha_0 = \alpha_1 = 0$ and $1 < p_0 < p_1 < \infty$ we obtain from (\ref{Cesthetap}) that
$$
\big ( Ces_{p_0}, Ces_{p_1} \big)_{\theta, p} = Ces_{p} ~ {\rm for} ~~ \frac{1}{p} = \frac{1-\theta}{p_0} + \frac{\theta}{p_1},
$$
which was already proved in \cite[Corollary 3.2]{AM13} and \cite[Theorem 2]{AM14a} using the identification from (\ref{CesaroK}) 
$(L^1, L^1(1/s))_{1-1/p, p} = Ces_p \, (1 < p < \infty)$ and reiteration theorem for the $K$-method $(\cdot)_{\theta, p}$ (see also 
\cite[Corollary 2]{Si91}).

From Theorem \ref{thm:tyldaspace} we obtain that if $X_0, X_1$ are Banach ideal spaces with the Fatou property such that either
$X_0, X_1$ are symmetric spaces or $C, C^*$ are bounded on both  $X_0$ and $X_1$, then 
$\widetilde{X_0} + \widetilde{X_1} = \widetilde{X_0 + X_1}$, which gives that 
$K(t, f; \widetilde{X_0}, \widetilde{X_1}) \approx K(t, \widetilde{f}; X_0, X_1)$ and we get commutativity of Tandori spaces 
with the real method of interpolation
\begin{equation} \label{4.5}
(\widetilde{X_0}, \widetilde{X_1})_G^K = [(X_0, X_1)_G^K]^{\Large  \sim}.
\end{equation}

%%%%%%%%%%%%%%%%%%%%%%%%%%%% Section 5
\section{Additional remarks}

A fundamental problem in interpolation theory is the description of all interpolation spaces with respect to a given Banach pair. 
In particular, it is not so rare that for a given Banach couple, all interpolation spaces may be generated by K-method and such 
couples are referred to be {\it Calder\'on couples} or {\it Calder\'on pairs}. 

A Banach couple $\bar{X} = (X_0, X_1)$ is called a {\it Calder\'on couple} if every interpolation space between $\bar{X}$  is described 
by the $K$-method. Equivalently, if for each $f, g \in X_0 + X_1$ satisfying
$$
K(t, f; X_0, X_1) \leq K(t, g; X_0, X_1) ~ {\rm for ~all} ~~t > 0,
$$
there is a bounded operator $T: X_0+X_1 \rightarrow X_0+X_1$ such that $\max (\| T\|_{X_0 \rightarrow X_0}, 
\| T\|_{X_1 \rightarrow X_1}) < \infty$ and $Tg = f$ (cf. \cite{BK91}, Theorem 4.4.5).

There are many examples of Calder\'on couples and couples which are not Calder\'on. What can we say about this problem for 
Ces\`aro, Copson and Tandori spaces?

By the Cwikel's result (cf. \cite[Theorem 1]{Cw81}) we obtain that  $(Ces_p(I), Ces_q(I))$ is a Calder\'on couple for $1 < p < q < \infty$, since
$$
(L^1[0, \infty), Ces_{\infty}[0, \infty))_{1-1/p, p} = Ces_p[0, \infty)
$$ 
and 
$$
(L^1(1-x)[0, 1], Ces_{\infty}[0, 1])_{1-1/p, p} = Ces_p[0, 1]
$$
(see Astashkin-Maligranda \cite[Proposition 3.1 and Theorem 4.1]{AM13}).

Moreover, using the Brudny{\v \i}-Dmitriev-Ovchinnikov theorem (see \cite[Theorem 1]{DO79} and \cite[Theorem 4.3.1]{BK91}) 
and Theorem 5 we obtain that if $X_0, X_1$ are Banach function spaces on $I = [0, \infty)$ such that both $C$ and $C^*$ are bounded on 
$X_0$ and on $X_1$, then $(CX_0, CX_1)$ is a Calder\'on couple.

Let us notice also that Masty{\l}o-Sinnamon \cite{MS06} proved that $(L^1, Ces_{\infty})$ is a Calder\'on couple and Le\'snik \cite{Le15} 
showed that also $(\widetilde{L^1}, L^{\infty})$ is a Calder\'on couple.
\vspace{2mm}

The following problems are natural to formulate here.
\vspace{2mm}

{\bf Problem 1}. For $I = [0, 1]$ and $1 \leq p < \infty$ identify $(L^1)^{1-\theta}(Ces_p)^{\theta}$ or $\varphi(L^1, Ces_p)$, 
or even more generally, $\varphi(L^1, CX)$.
\vspace{2mm}

{\bf Problem 2}. For $I = [0, 1]$ identify $(Ces_1)^{1-\theta}(Ces_{\infty})^{\theta}$ or $\varphi (Ces_1, Ces_{\infty})$, or even more general $\varphi(Ces_1, CX)$.

Note that  $(Ces_1)^{1-1/p}(Ces_{\infty})^{1/p} \neq Ces_p$ for $1 < p < \infty$. In fact, from (\ref{3.12}) we have that 
$(L^1)^{1-1/p}(Ces_{\infty})^{1/p} = Ces_p$ and by the uniqueness theorem (cf. \cite[Theorem 3.5]{CN03} or \cite[Corollary 1]{BM05})
$$
Ces_p = (L^1)^{1-1/p}(Ces_{\infty})^{1/p} \neq (Ces_1)^{1-1/p}(Ces_{\infty})^{1/p},
$$
since $L^1 \neq Ces_1$. Under some mild conditions on $\varphi$, from the uniqueness theorem proved in \cite[Theorem 1]{BM05}, we 
get that even $C[\varphi (L^1, L^{\infty})] \neq \varphi( Ces_1, Ces_{\infty})$.
\vspace{2mm}

{\bf Problem 3}. For $I = [0, 1]$ or $I = [0, \infty)$ identify $\widetilde{X}^{1-\theta} X^{\theta}$ or $\varphi(\widetilde{X}, X)$. 
\vspace{2mm}

The last identification can be useful for factorization since there appeared for $1 < p < \infty$ the $p$-convexification of  $Ces_{\infty}$, 
that is, $Ces_{\infty}^{(p)}$ and we have equalities (cf. \cite{KLM14})
$$
Ces_{\infty}^{(p)} = Ces_{\infty}^{1/p} (L^{\infty})^{1-1/p}= \big[ \widetilde{L^1}^{1/p} (L^{1})^{1-1/p} \big]^{\prime}.
$$  
Eventual identification in Problem 3 will suggests how to generalize factorization results presented in \cite{AM09} and \cite{KLM14}.
\vspace{2mm}

{\bf Problem 4}. What is an analogue of Theorem \ref{thm:Cesaro=real} for $I=[0,1]$?

%%%%%%%%%%%%%%%%%%%%%%%%%%%%%%%%%%%%%%%%%%%


\begin{thebibliography}{99}

\bibitem[Al57]{Al57} A. Alexiewicz, {\it On Cauchy's condensation theorem}, Studia Math. 16 (1957), 80--85. 

\bibitem[AM09]{AM09} S.V. Astashkin and L. Maligranda, {\it Structure of Ces\'aro function spaces}, Indag. Math. (N.S.) 20 
(2009), no. 3, 329--379.

\bibitem[AM13]{AM13} S.V. Astashkin and L. Maligranda, {\it Interpolation of Ces{\`a}ro sequence and function spaces}, 
Studia Math. 215 (2013), no. 1,  39--69.

\bibitem[AM14a]{AM14a} S.V. Astashkin and L. Maligranda, {\it Interpolation of Ces{\`a}ro and Copson spaces}, in: ``Banach and Function Spaces IV", Proc. of the Fourth Internat. Symp. on Banach and Function Spaces (ISBFS2012) (12-15 Sept. 2009, Kitakyushu-Japan), Edited by M. Kato, L. Maligranda and T. Suzuki, Yokohama Publishers 2014, 123--133.

\bibitem[AM14b]{AM14b} S.V. Astashkin and L. Maligranda, {\it Structure of Ces{\`a}ro function spaces: a survey}, Banach Center Publ. 
102 (2014), 13--40.

\bibitem[BS88]{BS88} C. Bennett and R. Sharpley, {\it Interpolation of Operators,} Academic Press, Boston 1988.

\bibitem[Be96]{Be96} G. Bennett, {\it Factorizing Classical Inequalities}, Mem. Amer. Math. Soc. 120 (1996), no. 576, 138 pp. 

\bibitem[Be81]{Be81} E. I. Berezhnoi, {\it Interpolation of positive operators in the spaces $\varphi (X_0, X_1)$}, In: Qualitative and Approximate Methods for the Investigation of Operator Equations. Yaroslav Gos. Univ., Yaroslavl 1981, 3--12 (in Russian).

\bibitem[BM05]{BM05} E. I. Berezhnoi and L. Maligranda, {\it Representation of Banach ideal spaces and factorization of operators}, 
Canad. J. Math. 57 (2005), 897--940.

\bibitem[BL76]{BL76} J. Bergh and J. L\"ofstr\"om, {\it Interpolation Spaces}, Springer, Berlin 1976. 

\bibitem[BK91]{BK91} Yu. A. Brudny{\u \i} and N. Ya. Krugljak, {\it Interpolation Functors and Interpolation Spaces}, 
North-Holland, Amsterdam 1991.

\bibitem[Ca64]{Ca64} A. P. Calder\'{o}n, {\it Intermediate spaces and interpolation, the complex method}, Studia Math. 
24 (1964), 113--190.

\bibitem[CR13]{CR13} G. P. Curbera and W. J. Ricker, {\it A feature of averaging}, Integral Equations Operator Theory 76 
(2013), no. 3, 447--449.

\bibitem[Cw81]{Cw81} M. Cwikel, {\it Monotonicity properties of interpolation spaces. II}, Ark. Mat. 19 (1981), 123--136.

\bibitem[Cw10]{Cw10} M. Cwikel, {\it Complex interpolation of compact operators mapping into lattice couples}, Proc. Est. 
Acad. Sci. 59 (2010),  no. 1, 19--28.

\bibitem[CN03]{CN03} M. Cwikel and P. Nilsson, {\it Interpolation of weighted Banach lattices},
Technion Preprint Series 834, Haifa 1989; Published as M. Cwikel, P. G. Nilsson and G. 
Schechtman, {\it Interpolation of Weighted Banach Lattices. A Characterization of Relatively 
Decomposable Banach Lattices}, Mem. Amer. Math. Soc. 165 (2003), no. 787,1--127. 

\bibitem[DS07]{DS07} O. Delgado and J. Soria, {\it Optimal domain for the Hardy operator}, J. Funct. Anal. 244 (2007), 
no. 1, 119--133.

\bibitem[DO79]{DO79} V. I. Dmitriev and V. I. Ovchinnikov, {\it Interpolation in spaces of the real method}, Dokl. Akad. 
Nauk SSSR 246 (1979), no. 4, 794--797; English transl. in: Soviet Math. Dokl. 20 (1979), no. 3, 538--542.

\bibitem[HLP52]{HLP52} G. H. Hardy, J. E. Littlewood and G. P\'olya, {\it Inequalities}, Cambridge Univ. Press 1952.

\bibitem[KA77]{KA77} L. V. Kantorovich and G. P. Akilov, {\it Functional Analysis}, Nauka, Moscow 1977 (Russian); English 
transl. Pergamon Press, Oxford-Elmsford, New York 1982. 

\bibitem[KMS07]{KMS07} R. Kerman, M. Milman and G. Sinnamon, {\it On the Brudny{\u\i}-Krugljak duality theory of spaces 
formed by the $K$-method of interpolation}, Rev. Mat. Complut. 20 (2007), no. 2, 367---389.

\bibitem[KL10]{KL10} P. Kolwicz and K. Le\'{s}nik, {\it Topological and geometrical structure of Calder\'on-Lozanovski{\u\i} 
construction}, Math. Inequal. Appl. 13 (2010), no. 1, 175--196. 

\bibitem[KLM13]{KLM13} P. Kolwicz, K. Le\'{s}nik and L. Maligranda, {\it Pointwise multipliers of Calder\'{o}n-Lozanovski{\u \i} 
spaces}, Math. Nachr. 286 (2013), no. 8-9, 876--907.

\bibitem[KLM14]{KLM14} P. Kolwicz, K. Le\'snik and L. Maligranda, {\it Pointwise products of some Banach function spaces 
and factorization}, J. Funct. Anal. 266 (2014), no. 2, 616--659.

\bibitem[KPS82]{KPS82} S. G. Krein, Yu. I. Petunin, and E. M. Semenov, {\it Interpolation of Linear Operators}, Amer. Math. 
Soc., Providence 1982.

 \bibitem[KMP93]{KMP93} N. Ya. Kruglyak, L. Maligranda and L. E. Persson, {\it A Carlson type inequality with blocks and 
 interpolation}, Studia Math. 104(1993), 161--180. 

\bibitem[KMP07]{KMP07} A. Kufner, L. Maligranda and L. E. Persson, {\it The Hardy inequality -- About its History and 
Some Related Results}, Vydavatelski Servis Publishing House, Pilzen 2007.

\bibitem[Le15]{Le15} K. Le\'snik, {\it Monotone substochastic operators and a new Calder\'on couple}, preprint of 14 pages, 
17 February 2015, {\tt arXiv:1502.04882} at: http://arxiv.org/pdf/1502.04882.pdf. 

\bibitem[LM15a]{LM15a} K. Le\'snik and L. Maligranda, {\it Abstract Ces\`aro spaces. Duality}, J. Math. Anal. Appl. 
424  (2015), no. 2, 932--951. 

\bibitem[LM15b]{LM15b} K. Le\'snik and L. Maligranda, {\it Abstract Ces\`aro spaces. Optimal range}, Integral Equations 
Operator Theory 81 (2015), no. 2, 227--235.

\bibitem[LT79]{LT79} J. Lindenstrauss and L. Tzafriri, {\it Classical Banach Spaces, II. Function Spaces}, Springer-Verlag, 
Berlin-New York 1979.

\bibitem[Lo69]{Lo69} G. Ya. Lozanovski\u{\i}, {\it On some Banach lattices}, Sibirsk. Mat. Zh. 10 (1969), no. 3, 584--599; 
English transl. in: Siberian. Math. J. 10 (1969), no. 3, 419--431.

\bibitem[Lo73]{Lo73} G. Ya. Lozanovski\u{\i}, {\it On some Banach lattices. IV}, Sibirsk. Mat. Zh. 14 (1973), 140--155; English 
transl. in: Siberian. Math. J. 14 (1973), 97--108.

\bibitem[Lo78a]{Lo78a} G. Ya. Lozanovski\u{\i}, \textit{Mappings of Banach lattices of measurable functions.} Izv. Vyssh. Uchebn. 
Zaved. Mat. 192 (1978), no. 5, 84--86; English transl. in: Soviet Math. (Iz. VUZ) 22 (1978), 61--63.

\bibitem[Lo78b]{Lo78b} G. Ya. Lozanovski\u{\i}, \textit{Transformations of ideal Banach spaces by means of concave functions}, in: 
Qualitative and Approximate Methods for the Interpolation of Operator Equations, No. 3, Yaroslav. Gos. Univ., Yaroslavl 1978, 
122--148 (in Russian).

\bibitem[LZ66]{LZ66} W. A. J. Luxemburg and A. C. Zaanen, {\it Some examples of normed K\"othe spaces}, Math. Ann. 
162 (1966), 337--350.

\bibitem[Ma85]{Ma85} L. Maligranda, {\it Calder\'on-Lozanovski{\u\i} spaces and interpolation of operators}, Semesterbericht Funktionalanalysis, T\"ubingen 8 (1985), 83--92.

\bibitem[Ma89]{Ma89} L. Maligranda, {\it Orlicz Spaces and Interpolation}, Seminars in Mathematics 5, University of Campinas, 
Campinas SP, Brazil 1989. 

\bibitem[MPS07]{MPS07} L. Maligranda, N. Petrot and S. Suantai, {\it On the James constant and $B$-convexity of 
Ces{\`a}ro and Ces{\`a}ro-Orlicz sequence spaces}, J. Math. Anal. Appl. {\bf 326} (2007), no. 1, 
312--331.

\bibitem[MS06]{MS06} M. Masty{\l}o and G. Sinnamon, {\it A Calder\'on couple of down spaces}, J. Funct. Anal. 240 (2006),  no. 1, 
192--225.

\bibitem[Ni85]{Ni85} P. Nilsson, {\it Interpolation of Banach lattices}, Studia Math. 82 (1985), 135--154.

\bibitem[Ov76]{Ov76} V. I. Ovchinnikov, {\it Interpolation theorems resulting from Grothendieck's inequality}, Funkcional. 
Anal. i Prilo\v zen. 10 (1976), no. 4, 45--54; English transl.: Functional Anal. Appl. 10 (1976), 287--294 (1977).

\bibitem[Ov84]{Ov84} V. I. Ovchinnikov, {\it The Method of Orbits in Interpolation Theory}, Math. Rep. 1 (1984),
no. 2, pp. i-x and 349--515.

\bibitem[RT10]{RT10} Y. Raynaud, and P. Tradacete, {\it Interpolation of Banach lattices and factorization of p-convex and 
q-concave operators}, Integral Equations Operator Theory 66 (2010), no. 1, 79--112. 

\bibitem[Re88]{Re88} S. Reisner, {\it On two theorems of Lozanovski{\u\i} concerning intermediate Banach lattices}, Lecture 
Notes in Math. 1317 (1988), 67--83. 

\bibitem[Ru80]{Ru80} Ja. B. Ruticki{\u \i}, {\it Operators with homogeneous kernels}, Sibirsk. Mat. Zh. 21 (1980), no. 1, 153--160; 
English transl. in: Siberian Math. J. 21 (1980), no. 1, 113--118.

\bibitem[Sh74]{Sh74} V. A. Shestakov, {\it On complex interpolation of Banach spaces of measurable functions}, Vestnik 
Leningrad. Univ. 19 (1974), 64--68; English transl. in Vestnik Leningrad Univ. Math. 7 (1979), 363--369 (1980). 

\bibitem[Sh81]{Sh81} V. A. Shestakov, {\it Transformations of Banach ideal spaces and interpolation of linear operators}, 
Bull. Polon. Acad. Sci. Math. 29 (1981), 569--577 (1982) (in Russian). 

\bibitem[Si91]{Si91} G. Sinnamon, {\it Interpolation of spaces defined by the level function}, in: Harmonic
Analysis (Sendai, 1990), ICM-90 Satell. Conf. Proc., Springer, Tokyo 1991, 190--193.

\bibitem[Ta55]{Ta55} K. Tandori, {\it {\"U}ber einen speziellen Banachschen Raum}, Publ. Math. Debrecen 3 (1954), 263--268 (1955).

\end{thebibliography}
\end{document}